\documentclass[reqno]{amsart}

\usepackage{amscd}

\theoremstyle{plain}
  \newtheorem{theorem}{Theorem}
  \newtheorem{corollary}{Corollary}
  \newtheorem{lemma}{Lemma}
  \newtheorem{proposition}{Proposition}
  \newtheorem*{conjecture*}{Conjecture}

\theoremstyle{definition}
  \newtheorem{definition}{Definition}
  
  \newtheorem*{convention*}{Convention}
  
\theoremstyle{remark}
  \newtheorem{remark}{Remark}
  \newtheorem{claim}{Claim}

\newcommand{\scalar}{\langle\cdot,\cdot\rangle}
\newcommand{\sscalar}{\langle\!\langle\cdot,\cdot\rangle\!\rangle}
\newcommand{\sm}{\text{\rm sm}}
\newcommand{\la}{\text{\rm la}}
\newcommand{\an}{\text{\rm an}}
\newcommand{\sing}{\text{\rm sing}}
\newcommand\R{\mathbb R}
\DeclareMathOperator{\ind}{ind}
\DeclareMathOperator{\GL}{GL}
\DeclareMathOperator{\ev}{ev}
\DeclareMathOperator{\Ev}{Ev}
\DeclareMathOperator{\length}{length}
\DeclareMathOperator{\dist}{dist}

\DeclareMathOperator{\Lef}{Lef}
\DeclareMathOperator{\sign}{sign}
\DeclareMathOperator{\Id}{Id}
\DeclareMathOperator{\pr}{pr}
\DeclareMathOperator{\IND}{IND}
\DeclareMathOperator{\Vol}{Vol}
\DeclareMathOperator{\Maps}{Maps}
\DeclareMathOperator{\Spect}{Spect}

\DeclareMathOperator{\Int}{Int}
\DeclareMathOperator{\grad}{grad}
\DeclareMathOperator{\Cr}{Cr}
\DeclareMathOperator{\tr}{tr}
\DeclareMathOperator{\diver}{div}
\DeclareMathOperator{\Gr}{Gr}
\DeclareMathOperator{\id}{id}
\DeclareMathOperator{\End}{End}
\DeclareMathOperator{\cs}{cs}

\def\cal{\mathcal}

\newcommand\itemref[1]{(\ref{#1})}

\begin{document}

\title[Laplace transform, Dynamics and Spectral geometry]
      {Laplace transform, Dynamics and Spectral geometry} 

\author{Dan Burghelea}

\address{Dan Burghelea,
         Dept. of Mathematics, 
         The Ohio State University, 
         231 West Avenue, Columbus, OH 43210, USA.}

\email{burghele@mps.ohio-state.edu}

\author{Stefan Haller}

\address{Stefan Haller,
         Department of Mathematics, University of Vienna,
         Nordbergstra{\ss}e 15, A-1090, Vienna, Austria.}

\email{stefan.haller@univie.ac.at}

\thanks{Part of this work was done while the second author enjoyed the
        warm hospitality of The Ohio State University. The second author is 
        supported by the \emph{Fonds zur F\"orderung der 
        wis\-sen\-schaft\-lichen Forschung} (Austrian Science Fund), 
        project number {\tt P14195-MAT}}




\begin{abstract} 
  We consider vector fields $X$ on a closed manifold $M$
  with rest points of Morse type. For such vector fields we define 
  the property of exponential growth. A cohomology class
  $\xi\in H^1(M;\mathbb R)$
  which is Lyapunov for $X$ defines counting functions 
  for isolated instantons and closed trajectories. If $X$ has  
  exponential growth property we show, under a mild hypothesis
  generically satisfied, that these counting functions can be recovered from 
  the spectral geometry associated to $(M,g,\omega)$ where $g$ is a Riemannian 
  metric and $\omega$ is a closed one form representing $\xi,$ cf Theorems 3 and 4
  in sectionn 1.6.
  This is done with the help of Dirichlet series and their Laplace transform.
\end{abstract}

\maketitle

\setcounter{tocdepth}{1}
\tableofcontents

\section{Introduction}
  \label{S:intro}

\subsection{Vector fields with zeros of Morse type and Lyapunov
            cohomology class}
\label{S:intro1}

Let $X$ be a smooth vector field on a smooth manifold $M$. A point
$x\in M$ is called a \emph{rest point} or a \emph{zero} of $X$ if $X(x)=0$. 
Denote by $\mathcal X:=\{x\in M|X(x)=0\}$ the set of rest points of $X.$

Recall that:
\begin{enumerate}
\item A \emph{parameterized trajectory} is a map $\theta :\mathbb R\to M$
so that $\theta'(t)= X(\theta (t)).$ A \emph{trajectory} 
is an equivalence 
class of parameterized trajectories with $\theta_1 \equiv \theta_2$
iff $\theta_1(t+ a)= \theta_2(t)$ for some real number $a$. Any representative 
$\theta$ of a trajectory
is called a parametrization.
\item
An \emph{instanton} from the rest point $x$ to the rest point $y$ is an isolated 
trajectory with the property that for one and then any parameterization $\theta,$ 
$\lim _{t\to-\infty}\theta(t)=x$, $\lim_{t\to+\infty}\theta(t)=y.$
\item 
A \emph{parameterized closed trajectory} is a pair $\tilde\theta=({\theta},T)$, 
with $\theta$ a parametrized trajectory and $T$ a positive real number 
so that $\theta(t+T)=\theta(t)$. 
A parameterized closed trajectory gives rise to a smooth map 
$\theta:S^1:=\mathbb R/T\mathbb Z\to M$. A \emph{closed trajectory}  
is an equivalence class $[\tilde \theta]$ of parameterized closed trajectories with 
$(\theta_1,T_1)\equiv(\theta_2,T_2)$
iff $\theta_1\equiv\theta_2$ and $T_1=T_2$.
\end{enumerate}

Recall that a rest point  $x\in\mathcal X$ is said to be of \emph{Morse type} if there
exist coordinates $(t_1,\dotsc,t_n)$ around $x$ so that 
$X=2\sum_{i=1}^qt_i\frac\partial{\partial t_i}-2\sum_{i=q+1}^nt_i\frac\partial{\partial t_i}.$
The integer $q$ is called the \emph{Morse index} of $x$ and denoted by  $\ind(x)$. A 
rest point of Morse type is non-degenerate and its Hopf index is $(-1)^{n-q}$.
It is independent of the chosen coordinates $(t_1,\dotsc,t_n).$ 
Then $\mathcal X=\bigsqcup_q\mathcal X_q$ where $\mathcal X_q$ denotes the set of rest points of index $q$.

For any rest point of Morse type $x,$ the \emph{stable} resp.\ \emph{unstable set} is defined by:
$$
W^\pm_x:=\{y|\lim_{t\to\pm\infty}\Psi_t(y)=x\}
$$
where $\Psi_t:M\to M$ denotes the flow of $X.$ The stable and unstable sets are
images of injective smooth immersions $i^\pm_x:W^\pm_x\to M$. By abuse of notation we denote 
the source manifold also by
$W^\pm_x$. The manifold $W^-_x$ resp.\ $W^+_x$ is diffeomorphic to $\mathbb R^{\ind(x)}$
resp.\ $\mathbb R^{n-\ind(x)}$.

\begin{convention*}
Unless explicitly mentioned all the vector fields in this paper are assumed 
to have all rest points of Morse type, hence isolated.
\end{convention*}

\begin{definition}\label{D:0} 
A vector field $X$ is said to have 
the \emph{exponential growth property at a rest point $x$} if for 
some (and then any) Riemannian metric $g$ there exists a 
positive constant $C$ so that $\Vol(D_r(x))\leq e^{Cr}$, 
for all $r\geq0$. Here $D_r(x)\subseteq W^-_x$ denotes the 
disk of radius $r$ with respect 
to the induced Riemannian metric $(i^-_x)^*g$ on $W_x^-$ centered at $x\in W^-_x.$ 
A vector field $X$ is said to have the \emph{exponential growth
property} if it has the exponential growth property at all
rest points.
\end{definition}

\begin{definition}\label{D:00}
A cohomology class $\xi\in H^1(M;\mathbb R)$ is called 
\emph{Lyapunov class} for a vector field $X$ if there exits a 
Riemannian metric $g$ and a closed one form $\omega$ 
representing $\xi$ so that $X=-\grad_g\omega$.
\end{definition}

\begin{remark}(to Definition 2)
1. An equivalent definition  is the following:
There exists a closed one form $\omega$ representing $\xi$
so that $\omega(X)<0$ on $M\setminus\mathcal X$ and such that
in a neighborhood of any rest point the vector field
$X$ is equal to $-\grad_g\omega$ for some Riemannian
metric $g$. It is proved in section~\ref{S:rho} that the two
definitions are actually equivalent.

2. The  closed form $\omega$ is a Morse form, i.e.\
locally it is the differential of a smooth function whose
critical points are non-degenerate.

3. Not all vector fields admit Lyapunov cohomology
classes. 
\end{remark}

\begin{definition}\label{D:MS}
The vector field $X$ is said to satisfy the \emph{Morse--Smale property}, MS for short,
if for any $x,y\in\mathcal X$ the maps $i^-_x$ and $i^+_y$ are transversal.
\end{definition}

We expect that every vector field which has a Lyapunov cohomology class,
and satisfies the Morse--Smale property, 
has the exponential growth property, cf.\ the conjecture
in section~\ref{SS:exp_rho}.
For the sake of Theorem 4 we introduce in section~\ref{S:homotopy}, 
cf.\ Definition~\ref{D:strong_exp},
the \emph{strong exponential growth property}. If the conjecture is true
both concepts are superfluous 
for the results of this paper.

In this paper we will show that a vector field $X$ and
a Lyapunov class $\xi$ for $X$
provide \emph{counting functions for the instantons} from $x$ to $y$ when 
$\ind(x)-\ind(y)=1$ 
and \emph{counting functions for closed trajectories}. Moreover these counting 
functions can be interpreted as Dirichlet series. 
If the vector field has exponential growth property 
these series have a finite abscissa of convergence, hence have a Laplace 
transform, cf section 1.2.
Their Laplace transform can be read off from the spectral geometry of a pair $(g,\omega)$ 
where $g$ is a Riemannian metric and $\omega$ is a closed one form representing $\xi$.

We will describe these counting functions and prove our results under the hypotheses that  
properties MS and NCT defined below are satisfied. Generically these properties 
are always satisfied, cf.\ Proposition~\ref{Prop:2} 
below.

Also in this paper, for any vector field $X$ and cohomology class $\xi\in H^1(M;\mathbb R)$ 
we define an invariant $\rho(\xi,X)\in \mathbb R\cup\{\pm\infty\}$ and show that if $\xi$ is 
Lyapunov for $X$ then 
exponential growth property is equivalent to $\rho(\xi,X)<\infty$.

If the vector field $X$ satisfies MS  then the set $\mathcal M(x,y)=W^-_x\cap W^+_y,$  
$x,y\in\mathcal X$ is 
the image by an injective immersion of a smooth manifold of dimension 
$\ind(x)-\ind(y)$ on which $\mathbb R$ acts freely. The quotient is a smooth manifold 
$\mathcal T(x,y)$ of dimension $\ind(x)-\ind(y)-1$
called the manifold of \emph{trajectories} from $x$ to $y$. If
$\ind(x)-\ind(y)=1$ then  $\mathcal T(x,y)$ is zero dimensional and its elements 
are isolated trajectories called \emph{instantons}. 

Choose $\mathcal O=\{\mathcal O_x\}_{x\in\mathcal X}$ a collection of orientations of the 
unstable manifolds of the critical points, with $\mathcal O_x$ an orientation of $W^-_x$. 
Any instanton $[\theta]$ from $x\in\mathcal X_q$ to $y\in\mathcal X_{q-1}$ has a sign
$\epsilon([\theta])=\epsilon^\cal O ([\theta])= \pm 1$ defined as follows: 
The orientations $\mathcal O_x$ and $\mathcal O_y$ induce an orientation
on $[\theta]$. Take $\epsilon([\theta])=+1$ if this orientation is compatible 
with the orientation from $x$ to $y$ and $\epsilon([\theta])=-1$ otherwise.

Let 
$\Psi_t$ denote the flow of $X.$ 
The closed trajectory $[\tilde \theta]$ is called non-degenerate if for some
(and then any) $t_0\in\mathbb R$ and representative $\tilde \theta= (\theta, T)$ the differential
$D_{\theta(t_0)}\Psi_T:T_{\theta(t_0)}M\to T_{\theta(t_0)}M$ is invertible with the
eigenvalue 1 of multiplicity one. 

\begin{definition}\label{D:NCT} 
The vector field $X$ is said to satisfies the \emph{non-degenerate closed trajectories
property}, NCT for short, if all closed trajectories of $X$ are
non-degenerate.
\end{definition}

Any non-degenerate closed trajectory $[\tilde\theta]$ has a \emph{period}
$p([\tilde\theta])\in\mathbb N$ and a \emph{sign} $\epsilon([\tilde\theta]):=\pm1$
defined as follows:
\begin{enumerate}
\item
$p([\tilde\theta])$ is the largest positive integer $p$ such that
$\theta:S^1\to M$ factors through a self map of $S^1$ of degree p.
\item
$\epsilon([\tilde\theta]):=\sign\det D_{\theta(t_0)}\Psi_T$ for some (and hence any) 
$t_0\in\mathbb R$ and parameterization $\tilde \theta.$
\end{enumerate}

A cohomology class $\xi\in H^1(M;\mathbb R)$ induces the homomorphism
$\xi:H_1(M;\mathbb Z)\to\mathbb R$ and then the injective group
homomorphism
$$
\text{$\xi:\Gamma_\xi\to\mathbb R$, with $\Gamma_\xi:=H_1(M;\mathbb Z)/\ker\xi$.}
$$

For any two points $x,y\in M$ denote by $\mathcal P_{x,y}$ the space of 
continuous paths from $x$ to $y$.
We say that $\alpha\in\mathcal P_{x,y}$ is equivalent to $\beta\in\cal P_{x,y}$,
iff the closed path $\beta^{-1}\star\alpha$ represents an element in $\ker\xi$.
(Here $\star$ denotes 
the juxtaposition of paths. Precisely if $\alpha,\beta:[0,1]\to M$ and
$\beta(0)=\alpha(1)$, then $\beta\star\alpha:[0,1]\to M$ is given by
$\alpha(2t)$ for $0\leq t\leq 1/2$ and $\beta(1-2t)$ for $1/2\leq t\leq 1$.)
 
We denote by $\hat{\mathcal P}_{x,y}=\hat{\mathcal P}^\xi_{x,y}$ the set 
of equivalence classes of elements in $\mathcal P_{x,y}$.
Note that $\Gamma_\xi$ acts freely and transitively,
both from the left and from the right, on $\hat{\mathcal P}^\xi_{x,y}$. 
The action $\star$ is defined by juxtaposing at $x$ resp.\ $y$ a closed 
curve representing an element $\gamma\in\Gamma_\xi$ to a path representing the 
element $\hat\alpha\in\hat{\mathcal P}_{x,y}^\xi$.

Any closed one form $\omega$ representing $\xi$ defines a map,
$\overline\omega:\mathcal P_{x,y}\to\mathbb R,$ by
$$
\overline\omega(\alpha):=\int_{[0,1]}\alpha^*\omega
$$
which in turn induces the map $\overline\omega:\hat{\mathcal P}^\xi_{x,y}\to\mathbb R$.
We have:
\begin {eqnarray*}
\overline\omega(\gamma\star\hat\alpha)
&=&
\xi(\gamma)+\overline\omega(\hat\alpha)
\\
\overline\omega(\hat\alpha\star\gamma)
&=& 
\overline\omega(\hat\alpha)+\xi(\gamma) 
\end{eqnarray*}
Note that for $\omega'=\omega+dh$ we have
$\overline{\omega'}=\overline{\omega}+h(y)-h(x)$.

\begin{proposition}\label{Prop:1}
Suppose $\xi\in H^1(M;\mathbb R)$ is a Lyapunov class for
the vector field $X$.
\begin{enumerate}
\item
If $X$ satisfies MS, $x\in\mathcal X_q$ and $y\in\mathcal X_{q-1}$ then the set of 
instantons from $x$ to $y$ in each class $\hat\alpha\in\hat{\mathcal P}^\xi_{x,y}$ is finite.
\item
If $X$ satisfies both MS and NCT then for any $\gamma\in\Gamma_\xi$
the set of closed trajectories representing the class $\gamma$ is finite.
\end{enumerate}
\end{proposition}

The proof is a straightforward consequence of the compacity of space of trajectories of bounded energy,
cf.\ \cite{BH01} and \cite{H02}.

Suppose $X$ is a vector field which satisfies MS and NCT and
suppose $\xi$ is a Lyapunov class for $X$.
In view of Proposition~\ref{Prop:1} we can define the 
\emph{counting function of closed trajectories} by
$$
\mathbb Z_X^\xi:\Gamma_\xi\to\mathbb Q,
\qquad
\mathbb Z_X^\xi(\gamma)
:=\sum_{[\tilde\theta]\in\gamma}
\frac{(-1)^{\epsilon([\tilde\theta])}}{p([\tilde\theta])}\in\mathbb Q.
$$  
If a collection of orientations $\mathcal O=\{\mathcal O_x\}_{x\in\mathcal  X}$ is given
one defines the \emph{counting function of the instantons} from $x$ to 
$y$ by
\begin{equation}\label{E:5}
\mathbb I^{X,\mathcal O,\xi}_{x,y}:\hat{\mathcal P}^\xi_{x,y}\to\mathbb Z,
\qquad
\mathbb I^{X,\mathcal O,\xi}_{x,y}(\hat\alpha)
:=\sum_{[\theta]\in\hat\alpha}\epsilon([\theta]).
\end{equation}
Note that the change of the orientations $\mathcal O$ might change the function 
$\mathbb I^{X,\mathcal O,\xi}_{x,y}$ but only up to multiplication by $\pm 1$. 
A key observation in this work is the fact that the counting functions 
$\mathbb I^{X,\mathcal O,\xi}_{x,y}$ and $\mathbb Z_X^\xi$ can be interpreted as Dirichlet series.

As long as Hypotheses MS and NCT are concerned we have the 
following genericity result. For a proof consult \cite{BH01} 
and the references in \cite[page 211]{Ha82}.

\begin{proposition}\label{Prop:2}
Suppose $X$ has $\xi\in H^1(M;\mathbb R)$ as a Lyapunov cohomology class.
\begin{enumerate}
\item
One can find a vector fields $X'$ arbitrarily close to $X$ in the $C^1$--topology which 
satisfy MS 
and have $\xi$ as Lyapunov cohomology class. Moreover one can choose $X'$ equal to $X$ 
in some neighborhood of $\cal X$ and away from any given neighborhood of $\cal X.$
\item
If in addition $X$ above satisfies MS one can find vector fields $X'$ 
arbitrary closed to $X$ in the $C^1$--topology which satisfy MS and NCT,  
and have $\xi$ as Lyapunov cohomology class.
Moreover one can choose $X'$ equal to $X$ 
in some neighborhood of $\cal X.$ 
\item
Consider the space of vector fields which have the same set of rest points as $X$,
and agree with $X$ in some neighborhood of $\cal X.$
Equip
this set with the $C^1$--topology. The subset of vector fields which satisfy MS and NCT 
is Baire residual set. 
\end{enumerate}
\end{proposition}

\subsection{Dirichlet series and their Laplace transform}
\label{S:intro6}

Recall that a Dirichlet series $f$ is given by a pair of finite or infinite sequences:
\begin{equation*}
\left(\begin{matrix}
\lambda_1 & < & \lambda_2 & < & \cdots & < & \lambda_k & < & \lambda_{k+1} & \cdots
\\
a_1       &   & a_2       &   & \cdots &   & a_k       &   & a_{k+1}       & \cdots
\end{matrix}\right)
\end{equation*}
The first sequence is a sequence of real numbers with the property that 
$\lambda_k\to\infty$ if the sequences are infinite. The second sequence is
a sequence of non-zero complex numbers. The associated series
\begin{equation*}
L(f)(z):=\sum_ie^{-z\lambda_i}a_i
\end{equation*}
has an \emph{abscissa of convergence}
$\rho(f)\leq \infty$,
characterized by the following properties, cf.\ \cite{Se79} and \cite{W46}:
\begin{enumerate}
\item
If $\Re z >\rho(f)$ then $f(z)$ is convergent and defines a 
holomorphic function.
\item
If $\Re z <\rho(f)$ then $f(z)$ is divergent.
\end{enumerate}
A Dirichlet series can be regarded as a complex valued measure with 
support on the discrete set 
$\{\lambda_1,\lambda_2,\dotsc\}\subseteq\mathbb R$
where the measure of $\lambda_i$ is equal to $a_i$. Then
the above series is the Laplace transform of this measure,
cf.\ \cite{W46}. The following proposition is a reformulation of results 
which lead to the Novikov theory and to the work of Hutchings--Lee and 
Pajitnov etc, cf.\ \cite{BH01} and \cite{H02} for more precise references.

\begin{proposition}\
\begin{enumerate}
\item
(Novikov)
Suppose $X$ is a vector field on a closed manifold $M$ 
which satisfies MS and has $\xi$ as a Lyapunov cohomology class.
Suppose $\omega$ is a closed one form representing $\xi$. 
Then for any $x\in\mathcal X_q$ and $y\in\mathcal X_{q-1}$
the collection of pairs of numbers
\begin{equation*}
\mathbb I^{X,\mathcal O,\omega}_{x,y}:=
\Bigl\{
\bigl(-\overline\omega(\hat\alpha),
\mathbb I^{X,\mathcal O,\xi}_{x,y}(\hat\alpha)\bigr)
\Bigm|
\mathbb I^{X,\mathcal O,\xi}_{x,y}(\hat\alpha)\neq 0, 
\hat\alpha\in\hat{\mathcal P}_{x,y}^\xi
\Bigr\}
\end{equation*}
defines a Dirichlet series. The sequence of $\lambda$'s
consists of the numbers $-\overline\omega(\hat\alpha)$
when $\mathbb I^{X,\mathcal O,\xi}_{x,y}(\hat\alpha)$ is non-zero, and the sequence
$a$'s consists of the numbers 
$\mathbb I^{X,\mathcal O,\xi}_{x,y}(\hat\alpha)\in\mathbb Z$.
\item
(D.~Fried, M.~Hutchings) 
If in addition $X$ satisfies NCT 
then the collection of pairs of numbers
\begin{equation*}
\mathbb Z^\xi_X:=
\Bigl\{
\bigl(-\xi(\gamma),\mathbb Z_X^\xi(\gamma)\bigr)
\Bigm|
\mathbb Z_X^\xi(\gamma)\neq0,\gamma\in\Gamma_\xi
\Bigr\}
\end{equation*}
defines a Dirichlet series. The sequence of $\lambda$'s
consists of the real numbers $-\xi(\gamma)$
when $\mathbb Z_X^\xi(\gamma)$ is non-zero and the sequence
of $a$'s consists of the numbers 
$\mathbb Z_X^\xi(\gamma)\in\mathbb Q$.
\end{enumerate}
\end{proposition}

We will show that if $X$ has exponential growth property then the Dirichlet series
$\mathbb I^{X,\mathcal O,\omega}_{x,y}$ and $\mathbb Z^\xi_X$ have the abscissa of 
convergence finite and therefore Laplace transform. 
The main results of this paper, Theorems~\ref{T:3} and \ref{T:4} below,
will provide explicit formulae for these Laplace transforms
in terms of the spectral geometry of $(M,g,\omega)$.
To explain such formulae we need additional considerations and results.

\subsection{The Witten--Laplacian}
\label{S:intro3}

Let $M$ be a closed manifold and $(g,\omega)$ a pair consisting of a 
Riemannian metric $g$ and a closed one form $\omega$. We suppose that 
$\omega$ is a \emph{Morse form}. This means that locally $\omega=dh$, $h$ 
smooth function with all critical points non-degenerate. A 
critical point or a zero of $\omega$ is a critical point of $h$
and since non-degenerate, has an index, the index of the Hessian $d^2_xh$, 
denoted by $\ind(x)$. Denote by $\mathcal X$
the set of critical points of $\omega$ and by $\mathcal X_q$ be the
subset of critical points of index $q$.

For $t\in\mathbb R$ consider the complex
$(\Omega^*(M),d^*_\omega(t))$ with differential
$d^q_\omega(t):\Omega^q(M)\to\Omega^{q+1}(M)$
given by
\begin{equation*}
d^q_\omega(t)(\alpha):=d\alpha+t\omega\wedge\alpha.
\end{equation*}
Using the Riemannian metric $g$ one constructs the formal adjoint of 
$d^q_\omega(t)$,
$
d^q_\omega(t)^\sharp:\Omega^{q+1}(M)\to\Omega^q(M),
$
and one defines the Witten--Laplacian
$\Delta^q_\omega(t):\Omega^q(M)\to\Omega^q(M)$ associated to the 
closed $1$--form $\omega$ by:
\begin{equation*}
\Delta^q_\omega(t):=d^q_\omega(t)^\sharp\circ d^q_t
+d^{q-1}_\omega(t)\circ d^{q-1}_\omega(t)^\sharp.
\end{equation*}
Thus, $\Delta^q_\omega(t)$ is a second order differential operator, with
$\Delta^q_\omega (0)=\Delta^q$, the Laplace--Beltrami operator. The operators
$\Delta^q_\omega(t)$ are elliptic, selfadjoint
and nonnegative, hence their spectra
$\Spect\Delta^q_\omega(t)$ lie in the interval $[0,\infty)$.
It is not hard to see that
\begin{equation*}
\Delta^q_\omega(t)=\Delta^q+t(L+L^\sharp)+t^2||\omega||^2\Id,
\end{equation*}
where $L$ denotes the Lie derivative along the vector field 
$-\grad_g\omega$, $L^\sharp$ the formal adjoint of $L$
and $||\omega||^2$ is the fiber wise norm of $\omega.$

The following result extends a result due to E.~Witten (cf.\ \cite{Wi82}) 
in the case that $\omega$ is exact and its proof was sketched in \cite{BH01}.

\begin{theorem}\label{T:1}
Let $M$ be a closed manifold and $(g,\omega)$ be a pair as above. Then there exist constants 
$C_1,C_2,C_3,T>0$ so that for $t>T$ we have:
\begin{enumerate}
\item\label{T:1:i}
$\Spect\Delta^q_\omega(t) \cap [C_1e^{-C_2 t}, C_3t]=\emptyset$.
\item\label{T:1:ii}
$\sharp\bigl(\Spect\Delta^q_\omega(t)\cap[0, C_1e^{-C_2 t}]\bigr)
=\sharp\mathcal X_q$.
\item\label{T:1:iii}
$1\in(C_1 e^{-C_2t},C_3t)$.
\end{enumerate}
Here $\sharp A$ denotes cardinality of the set $A$.
\end{theorem}


Theorem~\ref{T:1} can be complemented with the following proposition, see
Lemma~1.3 in \cite{BF97}.

\begin{proposition}\label{P:1}
For all but finitely many $t$ the dimension of 
$\ker\Delta^q_\omega(t)$ is constant in $t$.
\end{proposition}
Denote by $\Omega^*_\sm(M)(t)$ the $\mathbb R$--linear span of the 
eigen forms which correspond to eigenvalues smaller than $1$ referred 
bellow as the small eigenvalues. Denote by $\Omega^*_\la(M)(t)$ 
the orthogonal complement of $\Omega^*_\sm(M)(t)$ which, by elliptic 
theory, is a closed subspace of $\Omega^*(M)$ with respect to 
$C^\infty$--topology, in fact with respect to any Sobolev topology.
The space $\Omega^*_\la(M)(t)$ is the closure of the span of the 
eigen forms which correspond to eigenvalues larger than one.
As an immediate consequence of Theorem~\ref{T:1} we have for $t>T$ :
\begin{equation}\label{E:6}
\bigl(\Omega^*(M), d_\omega(t)\bigr)
=\bigl(\Omega^*_\sm(M)(t),d_\omega(t)\bigr)
\oplus\bigl(\Omega^*_\la(M)(t),d_\omega(t)\bigr)
\end{equation}
With respect to this decomposition the Witten--Laplacian is diagonalized
\begin{equation}\label{E:7}
\Delta^q_\omega(t)
=\Delta^q_{\omega,\sm}(t)\oplus\Delta^q_{\omega,\la}(t).
\end{equation}
and by Theorem~\ref{T:1}\itemref{T:1:ii}, we have for $t>T$
\begin{equation*}
\dim\Omega^q_\sm(M)(t)=\sharp\mathcal X_q.
\end{equation*}

The cochain complex $(\Omega^*_\la(M)(t),d_\omega(t))$ is acyclic 
and in view of Theorem~\ref{T:1}\itemref{T:1:ii} of finite codimension in the 
elliptic complex $(\Omega^*(M), d_\omega(t))$. Therefore we can 
define the function
\begin{equation}\label{E:6x}
\log T_{\an,\la}(t)
=\log T_{\an,\la}^{\omega,g}(t)
:=\frac12\sum_q (-1)^{q+1}q\log\det\Delta^q_{\omega,\la}(t)
\end{equation}
where $\det\Delta^q_{\omega,\la}(t)$ is the zeta-regularized
product of all eigenvalues of $\Delta^q_{\omega,\la}(t)$ larger 
than one.
This quantity will be referred to as the \emph{large analytic torsion}.

\subsection{Canonical base of the small complex}
\label{S:intro2}

Let $M$ be a closed manifold and $(g,g',\omega)$ be a triple consisting 
of two Riemannian metrics $g$ and $g'$ and a Morse form $\omega$. The vector field 
$X=-\grad_{g'}\omega$ has $[\omega]$ as a Lyapunov cohomology class. 

Suppose that $X$ satisfies MS and has exponential growth.
Choose $\mathcal O=\{\mathcal O_x\}_{x\in\mathcal X}$ a collection of 
orientations of the unstable manifolds  
with $\mathcal O_x$ orientation of $W^-_x$. 
Let $h_x:W^-_x\to\mathbb R$ be the unique smooth map defined by 
$dh_x=(i^-_x)^*\omega$ and $h_x(x)=0$. Clearly $h_x\leq 0$.

In view of the exponential growth property, cf.\ section~\ref{S:rho},
there exists $T$ so that for $t>T$ the integral 
\begin{equation}\label{E:9}
\Int^q_{X,\omega,\mathcal O}(t)(a)(x)
:=\int_{W^-_x}e^{th_x}(i_x^-)^*a,
\quad a\in\Omega^q(M),
\end{equation}
is absolutely convergent, cf.\ section~\ref{S:thm23}, and defines a linear map:
$$
\Int^q_{X,\omega,\mathcal O}(t):\Omega^q(M)\to\Maps(\mathcal X_q,\mathbb R).
$$

\begin{theorem}\label{T:2}
Suppose $(g,g',\omega)$ is a triple as above with $X$ of exponential growth 
and satisfying MS. Equip $\Omega^*(M)$ with the scalar product induced by 
$g$ and $\Maps(\mathcal X_q,\mathbb R)$ with the unique scalar product which 
makes $E_x\in\Maps(\mathcal X_q,\mathbb R)$, the characteristic functions of 
$x\in\mathcal X_q$, an orthonormal base.

Then there exists $T$ 
so that for any $q$ and $t\geq T$ the linear 
map $\Int^q_{X,\omega,\mathcal O}(t)$ defined by~\eqref{E:9}, when restricted 
to $\Omega^q_\sm(M)(t)$, is an isomorphism and an $O(1/t)$ isometry. 
In particular $\Omega^q_\sm(M)(t)$ has a 
canonical base $\{E^{\mathcal O}_x(t)|x\in\mathcal X_q\}$ with 
$E^{\mathcal O}_x(t)=(\Int^q_{X,\omega,\mathcal O}(t))^{-1}(E_x).$
\end{theorem}

As a consequence we have 
\begin{equation}\label{E:8}
d^{q-1}_{\omega}(E^\mathcal O_y(t))
=:\sum_{x\in\mathcal X_q}I^{X,\mathcal O,\omega,g}_{x,y}(t)\cdot E^\mathcal O_x(t),
\end{equation}
where $I^{X,\mathcal O,\omega,g}_{x,y}:[T,\infty)\to \mathbb R$
are smooth, actually analytic functions, cf.\ Theorem~\ref{T:3} below.

In addition to the functions 
$I^{X,\mathcal O,\omega,g}_{x,y}(t)$ defined for 
$t\geq T$, cf.\ \eqref{E:8}, we consider also the function
\begin{equation}\label{E:9x}
\log\mathbb V(t)
=\log\mathbb V_{\omega,g,X}(t)
:=\sum_q(-1)^q\log\Vol\{E_x(t)|x\in\mathcal X_q\}.
\end{equation}
Observe that the change in the orientations ${\mathcal O}$ does not change
the right side of \eqref{E:9x}, so ${\mathcal O}$ does not appear in the 
notation $\mathbb V(t)$.

\subsection{A geometric invariant associated to $(X,\omega,g)$
and a smooth function associated with the triple $(g,g',\omega)$}
\label{S:intro4}

Recall that Mathai--Quillen \cite{MQ86} (cf.\ also~\cite{BZ92}) have 
introduced a differential form
$\Psi_g\in\Omega^{n-1}(TM\setminus M;\mathcal O_M)$ for any
Riemannian manifold $(M,g)$ of dimension $n$. Here $\mathcal O_M$ denotes the orientation
bundle of $M$ pulled back to $TM$.
For any closed one form $\omega$ on $M$ we 
consider the  
form $\omega\wedge X^*\Psi_g\in\Omega^n(M\setminus\mathcal X;\mathcal O_M)$. 
Here $X=-\grad_{g'}\omega$ is regarded as a map $X:M\setminus\mathcal X\to TM\setminus M$ 
and $M$ is identified with the image of the zero section of the tangent bundle.

The integral
$$ 
\int_{M\setminus\mathcal X}\omega\wedge X^*\Psi_g
$$
is in general divergent. However it
does have a regularization defined by the formula
\begin{equation}\label{E:10}
R(X,\omega,g):=
\int_M\omega_0\wedge X^*\Psi_g
-\int_MfE_g
+\sum_{x\in\mathcal X}(-1)^{\ind(x)}f(x)
\end{equation}
where 
\begin{enumerate}
\item 
$f$ is a smooth function whose differential $df$ is equal to $\omega$ 
in a small neighborhood of $\mathcal X$ and therefore 
$\omega_0:=\omega-df$ vanishes in a small neighborhood of $\mathcal X$ and
\item 
$E_g\in\Omega^n(M;\mathcal O_M)$ is the Euler form associated with $g$.
\end{enumerate}
It will be shown in section~\ref{S:reg} below that the definition is 
independent of the choice of $f$, see also \cite{BH03}. Finally we introduce 
the function 
\begin {equation}\label{E:14}
\log\hat T_\an^{X,\omega,g}(t)
:=\log T_{\an,\la}^{\omega,g}(t)
-\log\mathbb V_{\omega,g,X}(t)
+ tR(X,\omega,g)
\end{equation}
where $X=-\grad_{g'}\omega$. 

\subsection{The main results}
\label{S:intro5}

The main results of this paper are Theorems~\ref{T:3} and 
\ref{T:4} below.

\begin{theorem}\label{T:3}
Suppose $X$ is a vector field which is MS and has exponential 
growth and suppose $\xi$ is a Lyapunov cohomology class for $X$.
Let $(g,g',\omega)$ be a system as in Theorem~\ref{T:2} so 
that $X=-\grad_{g'}\omega$ and $\omega$ a Morse form 
representing $\xi$. Let 
$I^{X,\mathcal O,\omega,g}_{x,y}:[T,\infty)\to\mathbb R$ 
be the functions defined by~\eqref{E:8}.
Then the Dirichlet series 
$\mathbb I^{X,\mathcal O,\xi}_{x,y}$ 
have finite abscissa of convergence and their Laplace 
transform are exactly the functions 
$I^{X,\mathcal O,\omega,g}_{x,y}(t)$. In particular 
$I^{X,\mathcal O,\omega,g}_{x,y}(t)$ is the restriction 
of a holomorphic function on $\{z\in\mathbb C|\Re z > T\}$.
\end{theorem}

\begin{theorem}\label{T:4}
Suppose $X$ is a vector field with $\xi$ a 
Lyapunov cohomology class which satisfies MS and NCT. Let 
$(g,g',\omega)$ be a system as in Theorem~\ref{T:2} so that 
$X=-\grad_{g'}\omega$ and $\omega$ a Morse form representing $\xi$.
Let $\log\hat T_\an^{X,\omega,g}(t)$ be the function defined by~\eqref{E:14}.

If in addition $X$ has exponential growth and $H^*(M,t[\omega])=0$ for $t$ 
sufficiently large
or $X$ has strong exponential growth 
then the Dirichlet series $\mathbb Z_X$ has finite abscissa of convergence and its Laplace 
transform is exactly the function  $\log\hat T_\an^{X,\omega,g}(t)$. In particular 
$\log\hat T_\an^{X,\omega,g}(t)$ is the restriction 
of a holomorphic function on $\{z\in\mathbb C|\Re z>T\}$.
\end{theorem}

(One can replace the acyclicity hypothesis by "$H^*_\sing(M;\Lambda_{\xi,\rho})$ is a 
free $\Lambda_{\xi,\rho}-$module 
for $\rho$ large enough".) 

If the conjecture in section~\ref{SS:exp_rho} is true, then the additional 
hypothesis (exponential growth resp.\ strong exponential growth) are 
superfluous.

\begin{remark}
The Dirichlet series $\mathbb Z_X$ depends only on $X$ and $\xi=[\omega]$,
while $\mathbb I^{X,\mathcal O,\xi}_{x,y}$ depends 
only on $X$ and $\xi$ up to multiplication with a constant 
(with a real number $r$ for the sequence of $\lambda$'s
and with $\epsilon=\pm1$ for the sequence of $a$'s).
\end{remark}

\begin{corollary}[J.~Marcsik cf.\ \cite{M98} or \cite{BH03}]\label{C:1}
Suppose $X$ is a vector field with no rest points,
$\xi\in H^1(M;\mathbb R)$ a Lyapunov class for $X$, $\omega$ a closed one form
representing $\xi$ and let $g$ a Riemannian metric on $M$.
Suppose all closed trajectories of $X$ are non-degenerate and denote by
$$
\log T_\an(t):=1/2\sum(-1)^{q+1}q\log\det(\Delta_\omega^q(t)).
$$
Then
$$
\log T_\an(t) + t\int_M\omega\wedge X^*\Psi_g
$$
is the Laplace transform of the Dirichlet series $\mathbb Z_X$ which
counts the set of closed trajectories of $X$ with the help of $\xi$.
\end{corollary}

\begin{remark}
In case that $M$ is the mapping torus of a diffeomorphism
$\phi:N\to N$, $M=N_\phi$ whose periodic points are all non-degenerate, 
the Laplace transform of the Dirichlet 
series $\mathbb Z_X$ is the Lefschetz zeta function $\Lef(Z)$ of $\phi$,
with the variable $Z$ replaced by $e^{-z}$.
\end{remark}

Theorems~\ref{T:3}, \ref{T:4} and Corollary~\ref{C:1} can be routinely 
extended to the case of a compact manifolds with boundary.

In section~\ref{S:top} we discuss one of the main topological tools in 
this paper, the completion of the unstable sets and of the space of 
unparameterized trajectories,
cf.\ Theorem~\ref{T:5}. This theorem was also proved in \cite{BH01}. 
In this paper we provide a significant short cut in the proof and a 
slightly more general formulation. 

In section~\ref{S:rho} we define the invariant $\rho$ and discuss the 
relationship with the exponential growth property.
Additional results of independent interest pointing toward the truth 
of the conjecture in section~\ref{SS:exp_rho} are also proved. 
The results of this section are not needed for the proofs of 
Theorems~\ref{T:2}--\ref{T:4}.

The proof of Theorem~\ref{T:1} as stated is contained in \cite{BH01} and so
is the proof of Theorem~\ref{T:2} but in a slightly different formulation
and (apparently) less generality. For this reason
and for the sake of completeness we will review and complete the
arguments (with  proper references to \cite{BH01} when necessary)
in section~\ref{S:thm23}.
Section~\ref{S:thm23} contains 
the proof of Theorem~\ref{T:2} and \ref{T:3}. 
Section~\ref{S:reg} treats the numerical invariant $R(X,\omega,g)$. 
The proof of 
Theorem~\ref{T:4} is presented in section~\ref{S:thm4} and 
relies on some previous work of Hutchings--Lee, Pajitnov
\cite{HL99}, \cite{H02}, \cite{P02} and the work of Bismut--Zhang and 
Burghelea--Friedlander--Kappeler \cite{BZ92}.

\section{Topology of the space of trajectories and unstable sets}
  \label{S:top}

In this section we discuss the 
completion of the unstable manifolds and of the manifolds of trajectories 
to manifolds with corners, which is a key topological tool in this work. The main result,
Theorem 5 is of independent interest.

\begin{definition}\label{P_defi}
Suppose $\xi\in H^1(M;\mathbb R)$. We say a covering
$\pi:\tilde M\to M$ satisfies 
\emph{property P with respect to $\xi$} if $\tilde M$ is
connected and $\pi^*\xi=0$.
\end{definition}

Let $X$ be vector field on a closed 
manifold $M$ which has $\xi\in H^1(M;\mathbb R)$ as a 
Lyapunov cohomology class, see Definition~\ref{D:00}. 
Suppose that $X$ satisfies MS. Let $\pi:\tilde M\to M$ be a covering
satisfying property P with respect to $\xi$.
Since $\xi$ is Lyapunov there exists a closed one form $\omega$
representing $\xi$ and a Riemannian metric $g$ so that 
$X=-\grad_g\omega$. Since the covering has property P we find
$h:\tilde M\to\mathbb R$ with $\pi^*\omega=dh$.

Denote by $\tilde X$ the vector field $\tilde X:=\pi^*X$. 
We write $\tilde{\mathcal X}=\pi^{-1}(\mathcal X)$ and 
$\tilde{\mathcal X}_q=\pi^{-1}(\mathcal X_q)$.
Clearly $\Cr(h)=\pi^{-1}(\Cr(\omega))$ are the zeros of 
$\tilde X$.

Given $\tilde x\in\tilde{\mathcal X}$ 
let $i^+_{\tilde x}:W^+_{\tilde x}\to\tilde M$
and $i^-_{\tilde x}:W^-_{\tilde x}\to\tilde M$, denote the 
one to one immersions whose images define the stable and 
unstable sets of $\tilde x$ with respect to the vector field 
$\tilde X$. The maps $i^\pm_{\tilde x}$ are actually smooth 
embeddings because $\tilde X$ is gradient like for the 
function $h$, and the manifold topology on 
$W^\pm_{\tilde x}$ coincides with the 
topology induced from $\tilde M$.
For any $\tilde x$ with $\pi(\tilde x)=x$ one can 
canonically identify $W^\pm_{\tilde x}$ to $W^\pm_x$
and then we have $\pi\circ i^\pm_{\tilde x}= i^\pm_x$.

As the maps $i^-_{\tilde x}$ and $ i^+_{\tilde y}$ are 
transversal, $\mathcal M(\tilde x,\tilde y):=
W^-_{\tilde x}\cap W^+_{\tilde y}$ is a submanifold of 
$\tilde M$ of dimension $\ind(\tilde x)-\ind(\tilde y)$. 
The manifold
$\mathcal M(\tilde x,\tilde y)$ is equipped with the action
$\mu:\mathbb R\times\mathcal M(\tilde x,\tilde y)\to\mathcal M(\tilde x,\tilde y)$,
defined by the flow generated by $\tilde X$. 
If $\tilde x\neq\tilde y$
the action $\mu$ is free and we denote the quotient
$\mathcal M(\tilde x,\tilde y)/\mathbb R$ by 
$\mathcal T(\tilde x,\tilde y)$.
The quotient $\mathcal T(\tilde x,\tilde y)$ is a smooth 
manifold of dimension
$\ind(\tilde x)-\ind(\tilde y)-1$, possibly empty, which, 
in view of the fact that $\tilde X(h) =\omega(X)<0$ is
diffeomorphic to the submanifold
$h^{-1}(c)\cap\mathcal M(\tilde x,\tilde y)$, where $c$ is 
any regular value of $h$ with $h(\tilde x)>c>h(\tilde y)$.

Note that if $\ind(\tilde x)\leq\ind(\tilde y)$, and 
$\tilde x\neq\tilde y$, in view the transversality required 
by the Hypothesis MS, the manifolds
$\mathcal M(\tilde x,\tilde y)$ and 
$\mathcal T(\tilde x,\tilde y)$ are empty.
We make the following convention: 
$\mathcal T(\tilde x,\tilde x):=\emptyset$.
This is very convenient for now 
$\mathcal T(\tilde x,\tilde y)\neq\emptyset$ implies 
$\ind(\tilde x)>\ind(\tilde y)$ and in particular 
$\tilde x\neq\tilde y$.

An \emph{unparameterized broken trajectory} from 
$\tilde x\in\tilde{\mathcal X}$ 
to $\tilde y\in\tilde {\mathcal X}$, is an element of the set
$\mathcal B(\tilde x,\tilde y)
:=\bigcup_{k\geq0}\mathcal B(\tilde x,\tilde y)_k$, where 
\begin{equation}\label{Bk_defi}
\mathcal B(\tilde x,\tilde y)_k
:=\bigcup
\mathcal T(\tilde y_0,\tilde y_1)
\times\cdots\times\mathcal T(\tilde y_k,\tilde y_{k+1})
\end{equation}
and the union is over all (tuples of) critical points 
$\tilde y_i\in\tilde{\mathcal X}$ with $\tilde y_0=\tilde x$ 
and $\tilde y_{k+1}=\tilde y$.

For $\tilde x\in \tilde{\mathcal X}$ introduce the
\emph{completed unstable set} $\hat W^-_{\tilde x}:=
\bigcup_{k\geq0}(\hat W^-_{\tilde x})_k$, where
\begin{equation}\label{Wk_defi}
(\hat W^-_{\tilde x})_k
:=\bigcup
\mathcal T(\tilde y_0,\tilde y_1)\times\cdots\times
\mathcal T(\tilde y_{k-1},\tilde y_k)\times W^-_{\tilde y_k}
\end{equation}
and the union is over all (tuples of) critical points 
$\tilde y_i\in\tilde{\mathcal X}$ with $\tilde y_0=\tilde x$.

To study $\hat W^-_{\tilde x}$ we introduce the set
$\mathcal B(\tilde x;\lambda)$ of
\emph{unparameterized broken trajectories from 
$\tilde x\in\tilde{\mathcal X}$ to the 
level $\lambda\in\mathbb R$} as $\mathcal B(\tilde x;\lambda):=
\bigcup_{k\geq0}\mathcal B(\tilde x;\lambda)_k$ where
$$
\mathcal B(\tilde x;\lambda)_k
:=\bigcup
\mathcal T(\tilde y_0,\tilde y_1)\times\cdots\times
\mathcal T(\tilde y_{k-1},\tilde y_k)\times
(W^-_{\tilde y_k}\cap h^{-1}(\lambda))
$$
and the union is over all (tuples of) critical points
$\tilde y_i\in\tilde{\mathcal X}$ with $\tilde y_0=\tilde x$.
Clearly, if $\lambda>h(\tilde x)$ then
$\mathcal B(\tilde x;\lambda)=\emptyset$.

Since any broken trajectory of $\tilde X$ intersects each level 
of $h$ in at most one point one can view the set 
$\mathcal B(\tilde x,\tilde y)$ resp.\ 
$\mathcal B(\tilde x;\lambda)$ as a subset of 
$C^0\big([h(\tilde y),h(\tilde x)],\tilde M\big)$ resp.\
$C^0\big([\lambda,h(\tilde x)],\tilde M\big)$. One parameterizes 
the points of a broken trajectory by the value of the function
$h$ on these points. This leads to the following 
characterization (and implicitly to a canonical
parameterization) of an unparameterized broken trajectory.

\begin{remark}
Let $\tilde x,\tilde y\in\tilde{\mathcal X}$ and set 
$a:=h(\tilde y)$, $b:=h(\tilde x)$. The parameterization 
above defines a one to one correspondence between
$\mathcal B(\tilde x,\tilde y)$ and the set of continuous 
mappings $\gamma:[a,b]\to\tilde M$, which satisfy the 
following two properties:
\begin{enumerate}
\item\label{R:3i}
$h(\gamma(s))=a+b-s$, 
$\gamma(a)=\tilde x$ and $\gamma(b)=\tilde{y}$.
\item\label{R:3ii}
There exists a finite collection of real numbers
$a=s_0<s_1<\cdots<s_{r-1}<s_r=b$, so that 
$\gamma(s_i)\in\tilde{\mathcal X}$ and $\gamma$ restricted to 
$(s_i,s_{i+1})$ has derivative at any point in the interval 
$(s_i,s_{i+1})$, and the derivative satisfies
\begin{equation*}
\gamma'(s)=\frac{\tilde X}{-\tilde X\cdot h}\big(\gamma(s)\big).
\end{equation*}
\end{enumerate}
Similarly the elements of $\mathcal B(\tilde x;\lambda)$ 
correspond to continuous mappings 
$\gamma:[\lambda,b]\to\tilde M$, which satisfies
\itemref{R:3i} and \itemref{R:3ii} with $a$ replaced by 
$\lambda$ and the condition $\gamma(b)=\tilde y$ ignored.
\end{remark}

We have the following proposition, which can be found 
in \cite{BH01}.

\begin{proposition}\label{P:comp}
For any $\tilde x,\tilde y\in\tilde{\mathcal X}$ and
$\lambda\in\mathbb R$, the spaces 
\begin{enumerate}
\item\label{P:comp:i}
$\mathcal B(\tilde x,\tilde y)$ with the topology induced from
$C^0\big([h(\tilde y),h(\tilde x)],\tilde M)$, and
\item\label{P:comp:ii}
$\mathcal B(\tilde x;\lambda)$ with the topology induced from
$C^0\big([\lambda,h(\tilde x)],\tilde M\big)$
\end{enumerate}
are compact.
\end{proposition}

Let $\hat i_{\tilde x}^-:\hat W^-_{\tilde x}\to\tilde M$ denote
the map whose restriction to 
$\mathcal T(\tilde y_0,\tilde y_1)
\times\cdots\times\mathcal T(\tilde y_{k-1},\tilde y_k)
\times W^-_{\tilde y_k}$
is the composition of the projection on $W^-_{\tilde y_k}$ with
$i_{\tilde y_k}^-$. Moreover let
$\hat h_{\tilde x}
:=h^{\tilde x}\circ\hat i_{\tilde x}^-:
\hat W^-_{\tilde x}\to\mathbb R$, 
where $h^{\tilde x}= h-h(\tilde x)$.

Recall that an $n$--dimensional \emph{manifold with corners} 
$P$, is a paracompact Hausdorff space equipped with a
maximal smooth atlas with charts 
$\varphi:U\to\varphi(U)\subseteq\mathbb R^n_+$,
where $\mathbb R^n_+=\{(x_1,\dotsc,x_n)\mid x_i\geq 0\}$. 
The collection of points of $P$ which correspond by some 
(and hence every) chart to points in $\mathbb R^n$ with 
exactly $k$ coordinates equal to zero is a well defined subset
of $P$ called the \emph{$k$--corner of $P$} and it will be 
denoted by 
$P_k$. It has a structure of a smooth $(n-k)$--dimensional 
manifold. The union $\partial P=P_1\cup P_2\cup\cdots\cup P_n$
is a closed subset which is a topological manifold and 
$(P,\partial P)$ is a topological manifold with boundary 
$\partial P$.

The following theorem was proven in~\cite{BH01} for the case
that $\tilde M$ is the minimal covering which has property P.

\begin{theorem}\label{T:5}
Let $M$ be a closed manifold, $X$ a vector field which is MS
and suppose $\xi$ is a Lyapunov class for $X$.
Let $\pi:\tilde M \to M$ be a covering which satisfies 
property P with respect to $\xi$ and let 
$h:\tilde M\to\mathbb R$ be a smooth map as above. Then:
\begin{enumerate}
\item\label{T:5i}
For any two rest points
$\tilde x,\tilde y\in\tilde{\mathcal X}$ the
smooth manifold $\mathcal T(\tilde x,\tilde y)$ has
${\mathcal B}(\tilde x,\tilde y)$ as a canonical
compactification. Moreover there is a canonic way to equip 
${\mathcal B}(\tilde x,\tilde y)$ 
with the structure of a compact smooth manifold with corners, 
whose $k$--corner is $\mathcal B(\tilde x,\tilde y)_k$
from \eqref{Bk_defi}.
\item\label{T:5ii}
For any rest point $\tilde x\in\tilde{\mathcal X}$, 
the smooth manifold $W^-_{\tilde x}$ has
$\hat W^-_{\tilde x}$ as a canonical completion. Moreover
there is a canonic way to equip
$\hat W_{\tilde x}^-$ with the structure of a smooth
manifold with corners, whose $k$--corner coincides with
$(\hat W^-_{\tilde x})_k$ from~\eqref{Wk_defi}.
\item\label{T:5iii}
$\hat i^-_{\tilde x}:\hat W^-_{\tilde x}\to\tilde M$
is smooth and proper, for all $\tilde x\in\tilde{\mathcal X}$.
\item\label{T:5iv}
$\hat h_{\tilde x}:\hat W^-_{\tilde x}\to\mathbb R$ 
is smooth and proper, for all $\tilde x\in\tilde{\mathcal X}$.
\end{enumerate}
\end{theorem}

\begin{proof} 
In view of Lemma~\ref{lemm3} in section 3, the set of Lyapunov classes for 
$X$ is open in $H^1(M;\mathbb R)$. So we can find a closed one
form $\omega$ and a Riemannian metric $g$ such that
$X=-\grad_g\omega$ and such that $\omega$ has degree of 
rationality one. Consider the minimal covering on which
$\xi=[\omega]$ becomes exact. Since $\xi$ has degree of 
rationality one
the critical values of $h$ form a discrete set. Recall that the closed 
one form $\omega$ has degree of 
rationality $k$ if the image of $[\omega](\Gamma)\subset\mathbb R$ is a 
free abelian group of rank $k.$
In \cite[paragraphs 4.1--4.3]{BH01} one can find all details
of the proof of Theorem~\ref{T:5} for this special $\xi$ and
this special covering.

Note that as long as properties \itemref{T:5i} through 
\itemref{T:5iii} are concerned they clearly remain true
when we pass to the universal covering of $M$ which obviously 
has property P. One easily concludes that they also remain
true for every covering which has property P. So we have
checked \itemref{T:5i} through \itemref{T:5iii} in the
general situation.

Next observe that $\hat h_{\tilde x}^-
=h^{\tilde x}\circ\hat i_{\tilde x^-}$ is certainly smooth as
a composition of two smooth mappings. The properness of
$\hat h_{\tilde x}$ follows from 
Proposition~\ref{P:comp}\itemref{P:comp:ii}.
\end{proof}

It will be convenient to formulate Theorem~\ref{T:5} without 
any reference to the covering $\pi:\tilde M\to M$ or to 
lifts $\tilde x$ of rest points $x$.

Let $\xi\in H^1(M;\mathbb R)$ be a one dimensional 
cohomology class so that $\pi^*\xi=0$.
As in section~\ref{S:intro}  denote by $\mathcal P_{x,y}$ 
the set of continuous paths from $x$ to $y$ and by 
$\hat{\mathcal P}_{x,y}^{\tilde M}$ the equivalence classes 
of paths in $\mathcal P_{x,y}$ with respect to the following 
equivalence relation.

\begin{definition}
Two paths $\alpha,\beta\in\mathcal P_{x,y}$ are equivalent if 
for some (and then for any) lift $\tilde x$ of 
$x$ the lifts $\tilde\alpha$ and $\tilde\beta$ of $\alpha$ and 
$\beta$ originating from $\tilde x$ end up
in the same point $\tilde y$.
\end{definition}

The reader might note that the present situation is slightly 
more general than the one considered in introduction which 
correspond to the case the covering $\pi$ is the 
$\Gamma_\xi$--principal covering with $\Gamma_\xi$ induced 
from $\xi$ as described in section~\ref{S:intro}. 
For this covering we have $\hat{\mathcal P}_{x,y}^{\tilde M}=
\hat{\mathcal P}_{x,y}^\xi$.

Note that any two lifts $\tilde x,\tilde y\in\tilde M$ 
determine an element $\hat\alpha\in\hat{\mathcal P}_{x,y}^{\tilde M}$ 
and the set of trajectories from $\tilde x$ to
$\tilde y$ identifies to the set $\mathcal T(x,y,\hat\alpha)$ 
of trajectories of $X$ from $x$ to $y$ in the 
class $\hat\alpha$.

Theorem~\ref{T:5} can be reformulated in the following way:

\begin{theorem}[Reformulation of Theorem~\ref{T:5}]\label{T:6}
Let $M$ be a smooth manifold, $X$ a smooth vector field 
which is MS and suppose $\xi$ is a Lyapunov class for $X$.
Let $\tilde M$ be a covering of $M$ which has property P with
respect to $\xi$. Then:
\begin{enumerate}
\item\label{T:6i}
For any two rest points $x,y\in\mathcal X$ and every
$\hat\alpha\in\hat{\mathcal P}_{x,y}^{\tilde M}$ the set 
$\mathcal T(x,y,\hat\alpha)$ has the structure 
of a smooth manifold of dimension $\ind(x)-\ind(y)-1$ 
which admits a canonical compactification to a compact smooth
manifold with corners $\mathcal B(x,y,\hat\alpha)$.
Its $k$--corner is 
$$
\mathcal B(x, y,\hat\alpha)_k
=\bigcup
\mathcal T(y_0,y_1,\hat\alpha_0)
\times\cdots\times\mathcal T(y_k, y_{k+1},\hat\alpha_k)
$$
where the union is over all (tuples of) critical points
$y_i\in\mathcal X$ and 
$\hat\alpha_i\in\hat{\mathcal P}_{y_i,y_{i+1}}^{\tilde M}$ with
$y_0=x$, $y_{k+1}=y$ and
$\hat\alpha_0\star\cdots\star\hat\alpha_k=\hat\alpha$.
\item\label{T:6ii}
For any rest point $x\in\mathcal X$ the smooth manifold
$W^-_x$ has a canonical completion to a smooth manifold
with corners $\hat W_x^-$. Its $k$--corner is
$$
(W^-_x)_k
=\bigcup
  \mathcal T( y_0, y_1,\hat\alpha_0)\times\cdots\times
  \mathcal T(y_{k-1},y_k,\hat\alpha_{k-1})\times W_{y_k}^-
$$
where the union is over all (tuples of) critical points
$y_i\in\mathcal X$ and $\hat\alpha_i
\in\hat{\mathcal P}_{y_i,y_{i+1}}^{\tilde M}$ with $y_0=x$.
\item\label{T:6iii}
The mapping $\hat i^-_x:\hat W^-_x\to M$ which on
$(W_x^-)_k$ is given by the composition of the 
projection onto $W_{y_k}^-$ with
$i^-_{y_k}:W^-_{y_k}\to M$ is smooth, for all 
$x\in\mathcal X$.
\item\label{T:6iv}
Let $\omega$ be a closed one form representing $\xi$.
Then the mappings $\hat h_x:\hat W^-_x\to\mathbb R$ which on 
$(W^-_x)_k$ is given by the composition of the projection
onto $W_{y_k}^-$ with $h^\omega_{y_k}:W^-_{y_k}\to\mathbb R$
plus 
$\overline\omega(\hat\alpha_0\star\cdots\star\hat\alpha_{k-1})$
is smooth and proper, for all $x\in\mathcal X$.
\end{enumerate}
\end{theorem}

The above results can be easily extended to the case of 
compact manifolds with boundary.

\subsection{Appendix to section 2}

Given a compact smooth manifold $M$ with boundary 
$\partial M$ we will consider only \emph{admissible metrics},
i.e.\ Riemannian metrics $g$ which are \emph{product like} 
near the boundary. In this case denote by $g_0$ the 
induced metric on the boundary. 
This means that there exists a collar neighborhood 
$\varphi:\partial M\times[0,\epsilon)\to M$
with $\varphi$ equal to the identity when restricted to 
$\partial M\times\{0\}$ and $\varphi^*g=g_0+ds^2$.

\begin{convention*}
Unless explicitly mentioned in this paper all the vector fields on a compact manifold with boundary
are assumed to be tangent to the boundary and have 
rest points of Morse type.
\end{convention*}

\begin{definition}\label{D:grad_with_bound}
The vector field $X$ has $\xi\in H^1(M;\mathbb R)$ as 
Lyapunov cohomology class if the following conditions hold:
\begin{enumerate}
\item
There exists a closed one form representing $\xi$ and an 
admissible metric so that $X=-\grad_g\omega$. In particular 
$X_{\partial M}=-\grad_{g_0}\omega_{\partial M}$,
where $\omega_{\partial M}$ denotes the 
pullback of $\omega$ to $\partial M$.
\item
If we set 
\begin{eqnarray*}
\mathcal X''_-
&:=&
\bigl\{x\in \mathcal X\cap \partial M\bigm| \ind_{\partial M}(x)=\ind(x)\bigr\}
\\
\mathcal X''_+
&:=&
\bigl\{x\in\mathcal X\cap\partial M\bigm|\ind_{\partial M}(x)=\ind(x)-1\bigr\}
\end{eqnarray*}
then $\mathcal X''_-$ resp.\ $\mathcal X''_+$ lie in different
components $\partial M_-$ resp.\ $\partial M_+$ of $M$. 
\end{enumerate}
\end{definition}

This definition implies that $\mathcal X=\mathcal X'\sqcup \mathcal X''$,
where $\mathcal X'$ is the set of rest points inside $M$ and $\mathcal X''$ of 
the rest points on $\partial M$ which is the same as the set of rest 
points of $X_{\partial M}$.
For $x\in \mathcal X''$ 
denote by $i^-_x:W^-_x \to M$ the unstable manifold with 
respect to $X$ and by $j^-_x:W^-_{\partial M,x}\to\partial M$ 
the unstable manifold with respect to $X_{\partial M}$.

\begin{remark}
\
\begin{enumerate}
\item
If $x\in\mathcal X''_-$ then the unstable manifold of $x$ with 
respect to $X$ and $X_{\partial M}$ are the same. More 
precisely $i^-_x:W^-_x\to M$ identifies to 
$j^-_x:W^-_{\partial M,x}\to\partial M$ followed by the 
inclusion of $\partial M \subset M$.
\item
If $x\in\mathcal X''_+$ then
\begin{enumerate}
  \item
  $(W^-_x,W^-_{\partial M,x})$ is a smooth manifold with 
  boundary diffeomorphic to $(\mathbb R^k_+,\mathbb R^{k-1})$ 
  with $k=\ind(x)$; and
  \item
  $i^-_x:W^-_x\to M$ is transversal to the boundary of $M$ 
  and $i^-_x: (i^-_x)^{-1}(\partial M)\to\partial M$ can be 
  identified to $j^-_x:W^-_{\partial M,x}\to\partial M$.
\end{enumerate}
\end{enumerate}
\end{remark}

Theorems~\ref{T:5} and \ref{T:6} remain true as stated 
with the following specifications. 
Set $P_y^-:=W^-_y\setminus W^-_{\partial M,y}$ for
$y\in \mathcal X''_+$, and $P_y^-:=W^-_y$ for
$y\in\mathcal X'\sqcup\mathcal X''_-$.
For $x\in\mathcal X''_+$ the $k$--corner of $\hat W_x^-$ then is
$$
(\hat W_x^-)_k=
(\hat W_{\partial M,x}^-)_{k-1}
\cup
\bigcup
\mathcal T(y_0,y_1,\hat\alpha_0)\times\cdots\times
\mathcal T(y_{k-1},y_k,\hat\alpha_{k-1})\times P_{y_k}^-
$$
where the big union is over all (tuples of) 
$y_i\in\mathcal X$ and 
$\hat\alpha_i\in\hat{\mathcal P}_{y_i,y_{i+1}}$ with 
$y_0=x$.

%

\section{Exponential growth property and the invariant $\rho$}
  \label{S:rho}

In this section we introduce for a pair $(X,\xi)$ consisting of a vector 
field $X$ and a cohomology class $\xi \in H^1(M;\mathbb R)$ 
an invariant $\rho(\xi,X)\in\mathbb R\cup\{\pm\infty\}$. 
For the purpose of this paper we are interested in the case this invariant 
is smaller than $\infty$. One expects that this is always the case if $\xi$ 
is Lyapunov for $X$ at least in the case $X$ satisfies MS.
If $X$ has $\xi$ as a Lyapunov cohomology class we prove that the exponential
growth and $\rho<\infty$ are equivalent. The discussion of this section is 
not needed for the proofs of Theorem~\ref{T:2}, \ref{T:3} and \ref{T:4}.

Throughout this section $M$ will be a closed manifold and $X$ a vector 
field with Morse zeros.

\subsection{The invariant $\rho$}
\label{SS:rho}

For a critical point $x\in\mathcal X$, i.e.\ a zero of $X$, we let 
$i_x^-:W_x^-\to M$
denote the smooth immersion of the unstable manifold into $M$. If $M$ is equipped with
a Riemannian metric we get an induced Riemannian metric $g_x:=(i_x^-)^*g$ on 
$W_x^-$ thus a volume density $\mu(g_x)$ on $W_x^-$ and hence the
spaces $L^p(W^-_x)$, $p\geq 1$. Though the $L^p$--norm depends on the metric
$g$ the space $L^p(W_x^-)$ and its topology does not. Indeed for another
Riemannian metric $g'$ on $M$ we find a constant $C>0$ so that
$1/C\leq\frac{g'(X,Y)}{g(X,Y)}\leq C$ for all tangent vectors $X$ and
$Y$ which implies 
$1/C'\leq\frac{\mu(g'_x)}{\mu(g_x)}\leq C'$ for some
constant $C'>0$.

Given a closed $1$--form $\omega$ on $M$ we let $h_x^\omega$ denote the unique
smooth function on $W_x^-$ which satisfies $dh_x^\omega=(i_x^-)^*\omega$ and
$h_x^\omega(x)=0$. We are interested in the space of $1$--forms for which
$e^{h^\omega_x}\in L^1(W_x^-)$. This condition actually only depends on the
cohomology class of $\omega$. Indeed we have
$h_x^{\omega+df}=h_x^\omega+(i_x^-)^*f-f(x)$ and so
$|h_x^{\omega+df}-h_x^\omega|\leq C''$ and 
$e^{-C''}\leq e^{h_x^{\omega+df}}/e^{h_x^\omega}\leq e^{C''}$
for some constant $C''>0$. So we define
$$
R_x(X)
:=\bigl\{
[\omega]\in H^1(M)
\bigm| 
e^{h^\omega_x}\in L^1(W_x^-)
\bigr\}
$$
and set $R(X):=\bigcap_{x\in\Cr(X)}R_x(X)$. 
Let us also define
$$
\rho_x(\xi,X):=\inf\{t\in\R\mid t\xi\in R_x(X)\}\in\R\cup\{\pm\infty\}
$$
as well as:
$$
\rho(\xi,X):=\inf\{t\in\R\mid t\xi\in R(X)\}\in\R\cup\{\pm\infty\}.
$$
Observe that $\rho(\xi,X)=\max_{x\in\mathcal X}\rho_x(\xi,X)$.

\begin{lemma}
The sets $R_x(X)$ and $R(X)$ are convex. Particularly
$$
\rho\bigl(\lambda\xi_1+(1-\lambda)\xi_2,X\bigr)
\leq\max\{\rho(\xi_1,X),\rho(\xi_2,X)\}
$$
for all $0\leq\lambda\leq 1$.
\end{lemma}

\begin{proof}
Indeed let $[\omega_0]$, $[\omega_1]\in H^1(M)$,
$\lambda\in[0,1]$ and set
$\omega_\lambda:=\lambda\omega_1+(1-\lambda)\omega_0$. Then
$h^{\omega_\lambda}_x=\lambda h^{\omega_1}_x+(1-\lambda)h^{\omega_0}_x$.
For $\lambda\in(0,1)$ we set $p:=1/\lambda>1$ and $q:=1/(1-\lambda)$. Then
$1/p+1/q=1$ and by H\"older's inequality
\begin{eqnarray*}
||e^{h^{\omega_\lambda}_x}||_1
&=&
||e^{\lambda h^{\omega_1}_x}e^{(1-\lambda)h^{\omega_0}_x}||_1
\\&\leq&
||e^{\lambda h^{\omega_1}_x}||_p||e^{(1-\lambda)h^{\omega_0}_x}||_q
\\&=&
||e^{h^{\omega_1}_x}||_1^\lambda||e^{h^{\omega_0}_x}||_1^{1-\lambda}
\end{eqnarray*}
So, if $[\omega_0]$ and $[\omega_1]\in R_x(X)$ then 
$[\omega_\lambda]\in R_x(X)$, and thus $R_x(X)$ is convex.
As an intersection of convex sets $R(X)$ is convex too.
\end{proof}

Next we introduce:
$$
B_x(X)
:=\bigl\{[\omega]\in H^1(M)\bigm| e^{h^\omega_x}\in L^\infty(W_x^-)\bigr\}
$$
and set $B(X):=\bigcap_{x\in\Cr(X)}B_x(X)$.
Note if $\xi$ is a Lyapunov cohomology class for $X$ then $\xi\in B(X)$,
cf.\ Lemma~\ref{lemm3} below.

Most obviously we have:

\begin{lemma}\label{lemm2}
The sets $B_x(X)$ and $B(X)$ are convex cones. Moreover we have
$R_x(X)+B_x(X)\subseteq R_x(X)$ and 
$R(X)+B(X)\subseteq R(X)$.
\end{lemma}

Next define
$$
L(X):=\bigl\{\xi\in H^1(M)\bigm| \text{$\xi$ is Lyapunov class for $X$}\bigr\}
$$
Recall from Definition~\ref{D:00} that $\xi\in L(X)$ if there exists a closed one form $\omega$ 
representing $\xi$ and a Riemannian metric $g$ such that $X=-\grad_g\omega$.

\begin{lemma}
Let $M$ be a smooth manifold, $X$ a vector field,
$\omega$ a closed one form and $g$ a Riemannian metric. 
Suppose $U\subset M$ is an open set and
\begin{enumerate}
\item 
the vector fields $X$ and $-\grad_g\omega$ agree on $U$ and
\item 
$\omega(X)<0$ on a neighborhood of $M\setminus U$.
\end{enumerate}
Then there exists a Riemannian metric $g'$ so that:
\begin{enumerate}
\item 
$X=-\grad_{g'}\omega$
\item 
$g$ and $g'$ agree on $U$.
\end{enumerate}
\end{lemma}

\begin{proof}
Let $N$ be an open neighborhood of $M\setminus U$ so that $\omega(X)<0$ 
and therefore $X_x\ne 0$, $x\in N$. For $x\in N$ the tangent space 
$T_x$ decomposes as the direct sum $T_xM=V_x\oplus[X_x]$
where $[X_x]$ denotes the one dimensional vector space generated by $X_x$ 
and $V_x=\ker(\omega(x):T_xM\to\mathbb R)$.
Clearly on $U$ the function $-\omega(X)$ is the square of the length of 
$X_x$ with respect to the metric $g$ and $X_x$
is orthogonal to $V_x$ and on $N$ it is strictly negative.
Define a new Riemannian metric $g'$ on $M$ as follows: 
For $x\in U$ the scalar product in $T_xM$ is the same as the one defined by $g$. 
For $x\in N$ the scalar product on $T_xM$ agrees to the one defined by
$g$ but make $V_x$ and $[X_x]$ perpendicular and the length of $X_x$ equal to 
$\sqrt{-\omega(X)(x)}$. It is clear that the new metric is well 
defined and smooth.
\end{proof}

\begin{corollary}\label{C:lyapunov}
Let $X$ be a vector field on $M$ and let $\xi\in H^1(M)$. 
Then $\xi$ is Lyapunov for $X$ if and only if there is a
closed one form $\omega$ representing $\xi$ and a 
Riemannian metric $g$ such that the following hold:
\begin{enumerate}
\item\label{gli}
$\omega(X)<0$ on $M\setminus\mathcal X$.
\item\label{glii}
$X=-\grad_g\omega$ on a neighborhood of $\mathcal X$.
\end{enumerate}
\end{corollary}


\begin{lemma}\label{lemm3}
The set $L(X)\subseteq H^1(M)$ is open and contained in 
$B(X)$. 
Moreover $L(X)$ is a convex cone.
\end{lemma}

\begin{proof}
The subset $L(X)\subseteq H^1(M)$ is open, for we can change the cohomology
class $[\omega]$ by adding a form whose support is disjoint from 
$\mathcal X$ and hence not affecting condition in 
Corollary~\ref{C:lyapunov}\itemref{glii}. 
If the form we add is sufficiently small the condition 
in Corollary~\ref{C:lyapunov}\itemref{gli}
will still be satisfied.

We have $L(X)\subseteq B_x(X)$ for $X=-\grad_g\omega$ implies that $h_x^\omega$
attains its maximum at $x$ and is thus bounded from above.

Next note that both conditions \itemref{gli} and \itemref{glii}
in Corollary~\ref{C:lyapunov} are convex
and homogeneous conditions on $\omega$. Thus $L(X)$ is a convex cone.
\end{proof}

\begin{lemma}
Suppose $L(X)\cap R(X)\neq\emptyset$. Then every ray of $L(X)$, i.e.\ a half
line starting at the origin which is contained in $L(X)$, intersects
$R(X)$.
\end{lemma}

\begin{proof}
Pick $\xi\in L(X)\cap R(X)$. 
Since $\xi\in R(X)$, Lemma~\ref{lemm3} and Lemma~\ref{lemm2} imply:
\begin{equation}\label{lemm3eq}
\xi+L(X)\subseteq R(X)+L(X)\subseteq R(X)+B(X)\subseteq R(X)
\end{equation}
On the other hand $L(X)$ is open and $\xi\in L(X)$,
so every ray in $L(X)$ has to intersect $\xi+L(X)$. In view of
\eqref{lemm3eq} it has to intersect $R(X)$ too.
\end{proof}

\begin{corollary}\label{C:rho_infty}
Suppose $\xi_0$ and $\xi$ are Lyapunov for $X$. Then
$\rho(\xi_0,X)<\infty$ implies $\rho(\xi,X)<\infty$.
\end{corollary}

\subsection{Exponential growth versus $\rho$}
\label{SS:exp_rho}

Let $x\in\mathcal X$ be a zero of $X$, $W_x^-$ the unstable manifold, let $g$ be a 
Riemannian metric on $M$ and let 
$r:=\dist(x,\cdot):W_x^-\to[0,\infty)$ denote the distance to $x$ with respect 
to the induced metric $g_x=(i^-_x)^*g$ on $W_x^-$. Clearly $r(x)=0$. Moreover let
$B_s(x):=\{y\in W^-_x|r(y)\leq s\}$ denote the ball of radius $s$ centered at $x$.

Recall from Definition~\ref{D:0} that
$X$ has the exponential growth property at a zero $x$ if 
there exists a constant $C\geq0$ such that 
$\Vol(B_s(x))\leq e^{Cs}$ for all $s\geq0$. Clearly 
this does not depend on the Riemannian metric $g$ on 
$M$ even though the constant $C$ will depend on $g$.

\begin{lemma}\label{L:exp:1}
Suppose we have $C\geq0$ such that $\Vol(B_s(x))\leq e^{Cs}$ for all $s\geq0$. Then 
$e^{-(C+\epsilon)r}\in L^1(W_x^-)$ for every $\epsilon>0$.
\end{lemma}

\begin{proof}
We have
\begin{equation}\label{ee1}
\int_{W_x^-}e^{-(C+\epsilon)r}
=
\sum_{n=0}^\infty\int_{B_{n+1}(x)\setminus B_n(x)}e^{-(C+\epsilon)r}
\end{equation}
On $B_{n+1}(x)\setminus B_n(x)$ we have $e^{-(C+\epsilon)r}\leq e^{-(C+\epsilon)n}$ and thus
\begin{eqnarray*}
\int_{B_{n+1}(x)\setminus B_n(x)}e^{-(C+\epsilon)r}
&\leq&
\Vol(B_{n+1}(x))e^{-(C+\epsilon)n}
\\&\leq&
e^{C(n+1)}e^{-(C+\epsilon)n}
=e^Ce^{-\epsilon n}
\end{eqnarray*}
So \eqref{ee1} implies
$$
\int_{W_x^-}e^{-(C+\epsilon)r}
\leq
e^C\sum_{n=0}^\infty e^{-\epsilon n}
=e^C(1-e^{-\epsilon})^{-1}<\infty
$$
and thus $e^{-(C+\epsilon)r}\in L^1(W_x^-)$.
\end{proof}

\begin{lemma}\label{L:exp:2}
Suppose we have $C\geq0$ such that $e^{-Cr}\in L^1(W_x^-)$. 
Then there exists a constant $C_0$ such that
$\Vol(B_s(x))\leq C_0e^{Cs}$ for all $s\geq0$.
\end{lemma}

\begin{proof}
We start with the following estimate for $N\in\mathbb N$:
\begin{eqnarray*}
\lefteqn{\Vol(B_{N+1}(x))e^{-C(N+1)}=}
\\&=&
\sum_{n=0}^N\Vol(B_{n+1}(x))e^{-C(n+1)}-\Vol(B_n(x))e^{-Cn}
\\&\leq&
\sum_{n=0}^\infty\bigl(\Vol(B_{n+1}(x))-\Vol(B_n(x))\bigr)e^{-C(n+1)}
\\&=&
\sum_{n=0}^\infty\Vol\bigl(B_{n+1}(x)\setminus B_n(x)\bigr)e^{-C(n+1)}
\\&\leq&
\sum_{n=0}^\infty\int_{B_{n+1}(x)\setminus B_n(x)}e^{-Cr}
=\int_{W^-_x}e^{-Cr}
\end{eqnarray*}
Given $s\geq0$ we choose an integer $N$ with $N\leq s\leq N+1$. Then
$\Vol(B_s(x))e^{-Cs}\leq\Vol(B_{N+1}(x))e^{-CN}=e^C\Vol(B_{N+1}(x))e^{-C(N+1)}$. So the computation above shows
$$
\Vol(B_s(x))e^{-Cs}
\leq e^C\int_{W_x^-}e^{-Cr}=:C_0<\infty
$$
and thus $\Vol(B_s(x))\leq C_0e^{Cs}$ for all $s\geq0$.
\end{proof}

As immediate consequence of the two preceding lemmas we have

\begin{proposition}\label{P:expgrow_vs_L1}
Let $x$ be a zero of $X$. Then the following are equivalent:
\begin{enumerate}
\item
$X$ has the exponential growth property at $x$ with respect to one (and hence every) 
Riemannian metric on $M$.
\item
For one (and hence every) Riemannian metric on $M$ there exists a constant $C\geq0$
such that $e^{-Cr}\in L^1(W_x^-)$.
\end{enumerate}
\end{proposition}

Let $g$ be a Riemannian metric on $M$, $\omega$ a closed one form and
consider $X=-\grad_g\omega$. Assume $X$ has all zeros of Morse type and let
$x$ be one of them. Recall that we have a smooth function 
$h_x^\omega:W_x^-\to(-\infty,0]$ defined by $(i_x^-)^*\omega=dh^\omega_x$ 
and $h_x^\omega(x)=0$. The next two lemmas tell, that 
$-h_x^\omega:W_x^-\to[0,\infty)$ is comparable with
$r:W_x^-\to[0,\infty)$.

\begin{lemma}\label{L:exp:3}
In this situation there
exists a constant $C_{\omega,g}\geq0$ such that $r\leq1-C_{\omega,g}h_x^\omega$ on $W_x^-$.
\end{lemma}

\begin{proof}
The proof is exactly the same as the one in \cite[Lemma~3(2)]{BH01}. 
Note that the MS property is not used there.
\end{proof}

\begin{lemma}\label{L:exp:4}
In this situation there exists a constant $C_{\omega,g}'\geq0$ such that
$-h_x^\omega\leq C_{\omega,g}'r$.
\end{lemma}

\begin{proof}
Let $\gamma:[0,1]\to W_x^-$ be any path starting at $\gamma(0)=x$. 
For simplicity set $h:=h_x^\omega$. Since $h(x)=0$ we find
$$
|h(\gamma(1))|
=
\Bigl|\int_0^1(dh)(\gamma'(t))dt\Bigr|
\leq
||\omega||_\infty\int_0^1|\gamma'(t)|dt
=
||\omega||_\infty\length(\gamma)
$$
with $||\omega||_\infty$ the supremums norm of $\omega$.
We conclude 
$$
||\omega||_\infty r(\gamma(1))
=
||\omega||_\infty\dist(x,\gamma(1))
\geq
|h(\gamma(1))|
\geq
-h(\gamma(1))
$$ 
and thus $-h\leq C_{\omega,g}'r$ with $C'_{\omega,g}:=||\omega||_\infty$.
\end{proof}

Let us collect what we have found so far.

\begin{proposition}\label{P:expgrow}
Suppose $\xi$ is Lyapunov for $X$ and let $x$ be a zero of $X$.
Then the following are equivalent.
\begin{enumerate}
\item\label{P:expgrow:i}
$X$ has the exponential growth property at $x$ with respect to one (and hence every) 
Riemannian metric on $M$.
\item\label{P:expgrow:ii}
For one (and hence every) Riemannian metric on $M$ there exists a constant $C\geq0$ such that
$e^{-Cr}\in L^1(W_x^-)$.
\item\label{P:expgrow:iii}
$\rho_x(\xi,X)<\infty$.
\end{enumerate}
\end{proposition}

\begin{proof}
The equivalence of \itemref{P:expgrow:i} and \itemref{P:expgrow:ii} was established in
Proposition~\ref{P:expgrow_vs_L1} without the assumption that $\xi$ is Lyapunov for $X$.
The implication \itemref{P:expgrow:ii} $\Rightarrow$ \itemref{P:expgrow:iii} follows from
Lemma~\ref{L:exp:3};
the implication \itemref{P:expgrow:iii} $\Rightarrow$ \itemref{P:expgrow:ii} from
Lemma~\ref{L:exp:4}.
\end{proof}

Note that this again implies Corollary~\ref{C:rho_infty}.
We expect that these equivalent statements hold true, at least in the generic
situation. More precisely:

\begin{conjecture*}[Exponential growth]\label{conj}
Let $g$ be a Riemannian metric on a closed manifold $M$, $\omega$ a closed 
one form (and assume $X=-\grad_g\omega$ is Morse--Smale.) 
Let $x$ be a zero and let
$i_x^-:W_x^-\to M$ denote its unstable manifold. Let $\Vol(B_s(x))$ denote the 
volume of the ball $B_s(x)\subseteq W_x^-$ of radius $s$ centered at $x\in W_x^-$
with respect to the induced Riemannian metric $(i_x^-)^*g$ on $W_x^-$. Then there exists a constant
$C\geq0$ such that $\Vol(B_s(x))\leq e^{Cs}$ for all $s\geq0$.
\end{conjecture*}

\subsection{A criterion for exponential growth}
\label{SS:exp_crit}

The rest of the section is dedicated to a criterion which guarantees
that the exponential growth property, and hence $\rho<\infty$, 
holds in simple situations.

Suppose $x\in\mathcal X_q$. Let $B\subseteq W_x^-$ denote a 
small ball centered at $x$.
The submanifold $i_x^-\bigl(W_x^-\setminus B\bigr)\subseteq M$ 
gives rise to a
submanifold $\Gr(W_x^-\setminus B)\subseteq\Gr_q(TM)$ 
in the Grassmannified tangent
bundle, the space of $q$--dimensional subspaces in $TM$. 
For a critical point $y\in\mathcal X$ we define
$$
K_x(y):=\Gr_q(T_yW_y^-)\cap\overline{\Gr(W_x^-\setminus B)}
$$
where $\Gr_q(T_yW_y^-)\subseteq\Gr_q(T_yM)\subseteq\Gr_q(TM)$.
Note that $K_x(y)$ does not depend on the choice of $B$.

\begin{remark}\label{rem_on_K}
\
\begin{enumerate}
\item
Even though we removed a neighborhood of $x$ from the
unstable manifold
$W_x^-$ the set $K_x(x)$ need not be empty. However if
we did not remove
$B$ the set $K_x(x)$ would never be vacuous for trivial
reasons.
\item\label{rem_on_K:ii}
If $q=\ind(x)>\ind(y)$ we have $K_x(y)=\emptyset$, for
$\Gr_q(T_yW_y^-)=\emptyset$ since $q>\dim(T_yW_y^-)=\ind(y)$.
\item\label{rhotop}
If $\dim(M)=n=q=\ind(x)$ we always have $K_x(y)=\emptyset$ for
all $y\in\mathcal X$.
\end{enumerate}
\end{remark}

\begin{proposition}\label{rhoxprop}
Let $\xi$ be Lyapunov for $X$ and suppose
$K_x(y)=\emptyset$ for all $y\in\mathcal X$. Then $\rho_x(\xi,X)<\infty$.
\end{proposition}

We start with a little

\begin{lemma}\label{little_lemm}
Let $(V,g)$ be an Euclidean vector space and $V=V^+\oplus V^-$ an orthogonal
decomposition. For $\kappa\geq0$ consider the endomorphism
$A_\kappa:=\kappa\id\oplus-\id\in\End(V)$ and the function
$$
\delta^{A_\kappa}:\Gr_q(V)\to\mathbb R,\quad
\delta^{A_\kappa}(W):=\tr_{g|_W}(p^\perp_W\circ A_\kappa\circ i_W),
$$
where $i_W:W\to V$ denotes the inclusion and $p^\perp_W:V\to W$ the
orthogonal projection. Suppose we have a compact subset $K\subseteq\Gr_q(V)$
for which $\Gr_q(V^+)\cap K=\emptyset$. Then there exists $\kappa>0$ and
$\epsilon>0$ with $\delta^{A_\kappa}\leq-\epsilon$ on $K$.
\end{lemma}

\begin{proof}
Consider the case $\kappa=0$. Let $W\in\Gr_q(V)$ and choose a 
$g|_W$--orthonormal base $e_i=(e^+_i,e^-_i)\in V^+\oplus V^-$, 
$1\leq i\leq q$, of $W$. Then
$$
\delta^{A_0}(W)=\sum_{i=1}^qg(e_i,A_0e_i)
=-\sum_{i=1}^qg(e_i^-,e_i^-).
$$
So we see that $\delta^{A_0}\leq 0$ and
$\delta^{A_0}(W)=0$ iff $W\in\Gr_q(V^+)$. Thus $\delta^{A_0}|_K<0$.
Since $\delta^{A_\kappa}$ depends continuously on $\kappa$ and since $K$ is
compact we certainly find $\kappa>0$ and $\epsilon>0$ so that
$\delta^{A_\kappa}|_K\leq-\epsilon$.
\end{proof}

\begin{proof}[Proof of Proposition~\ref{rhoxprop}]
Let $S\subseteq W_x^-$ denote a small sphere centered at $x$. 
Let $\tilde X:=(i_x^-)^*X$ denote the restriction of $X$ to $W_x^-$ and 
let $\Phi_t$ denote the flow of $\tilde X$ at time $t$. Then
$$
\varphi:S\times[0,\infty)\to W^-_x,\quad
\varphi(x,t)=\varphi_t(x)=\Phi_t(x)
$$
parameterizes $W_x^-$ with a small neighborhood of $x$ removed.

Let $\kappa>0$. For every $y\in\mathcal X$ choose a chart
$u_y:U_y\to\mathbb R^n$ centered at $y$ so that
$$
X|_{U_y}=\kappa
\sum_{i\leq\ind(y)}u_y^i\frac\partial{\partial u_y^i}
-\sum_{i>\ind(y)}u_y^i\frac\partial{\partial u_y^i}.
$$
Let $g$ be a Riemannian metric on $M$ which restricts to
$\sum_idu_y^i\otimes du_y^i$ on $U_y$ and set $g_x:=(i_x^-)^*g$. Then
$$
\nabla X|_{U_y}=\kappa
\sum_{i\leq\ind(y)}du_y^i\otimes\frac\partial{\partial u_y^i}
-\sum_{i>\ind(y)}du_y^i\otimes\frac\partial{\partial u_y^i}.
$$
In view of our assumption $K_x(y)=\emptyset$ for all $y\in\mathcal X$ 
Lemma~\ref{little_lemm} permits us to choose $\kappa>0$ and $\epsilon>0$ so
that after possibly shrinking $U_y$ we have
\begin{equation}\label{divestia}
\diver_{g_x}(\tilde X)
=\tr_{g_x}(\nabla\tilde X)
\leq-\epsilon<0
\quad\text{on}\quad 
\varphi(S\times[0,\infty))\cap(i_x^-)^{-1}
\Bigl(\bigcup_{y\in\mathcal X}U_y\Bigr).
\end{equation}
Next choose a closed $1$--form $\omega$ so that $[\omega]=\xi$ and
$\omega(X)<0$ on $M\setminus\mathcal X$. Choose $\tau>0$ so that
\begin{equation}\label{divestib}
\tau\omega(X)+\ind(x)||\nabla X||_g\leq-\epsilon<0
\quad\text{on}\quad
M\setminus\bigcup_{y\in\mathcal X}U_y.
\end{equation}
Using $\tau\tilde X\cdot h^\omega_x\leq 0$ and
$$
\diver_{g_x}(\tilde X)
=\tr_{g_x}(\nabla\tilde X)
\leq\ind(x)||\nabla\tilde X||_{g_x}
\leq\ind(x)||\nabla X||_g
$$
\eqref{divestia} and \eqref{divestib} yield 
\begin{equation}\label{divesti}
\tau \tilde X\cdot h^\omega_x+\diver_{g_x}(\tilde X)\leq-\epsilon<0
\quad\text{on}\quad
\varphi(S\times[0,\infty)).
\end{equation}
Choose an orientation of $W_x^-$ and let $\mu$ denote the volume form
on $W_x^-$ induced by $g_x$. Consider the function
$$
\psi:[0,\infty)\to\mathbb R,\quad
\psi(t):=\int_{\varphi(S\times[0,t])}e^{\tau h^\omega_x}\mu\geq 0.
$$
For its first derivative we find
$$
\psi'(t)=\int_{\varphi_t(S)}e^{\tau h^\omega_x}i_{\tilde X}\mu>0
$$
and for the second derivative, using \eqref{divesti},
\begin{eqnarray*}
\psi''(t)
&=&
\int_{\varphi_t(S)}
\bigl(\tau\tilde X\cdot h^\omega_x+\diver_{g_x}(\tilde X)\bigr)
e^{\tau h^\omega_x}i_{\tilde X}\mu
\\
&\leq&
-\epsilon\int_{\varphi_t(S)}e^{\tau h^\omega_x}i_{\tilde X}\mu
=-\epsilon\psi'(t).
\end{eqnarray*}
So $(\ln\circ\psi')'(t)\leq-\epsilon$ hence 
$\psi'(t)\leq\psi'(0)e^{-\epsilon t}$
and integrating again we find
$$
\psi(t)
\leq\psi(0)+\psi'(0)(1-e^{-\epsilon t})/\epsilon
\leq\psi'(0)/\epsilon.
$$
So we have $e^{\tau h^\omega_x}\in L^1\bigl(\varphi(S\times[0,\infty)\bigr)$
and hence $e^{\tau h^\omega_x}\in L^1(W_x^-)$ too. We conclude
$\rho_x(\xi,X)\leq\tau<\infty$.
\end{proof}

\begin{remark}\label{Rema:6}
\
\begin{enumerate}
\item
In view of Remark~\ref{rem_on_K}\itemref{rhotop} Proposition~\ref{rhoxprop}
implies $\rho_x(\xi,X)<\infty$ whenever $\xi$ is Lyapunov for
$X$ and $\ind(x)=\dim(M)$. However there is a much easier argument for this special
case. Indeed, in this case $W_x^-$ is an open subset of $M$ and therefore its
volume has to be finite. Since $\xi$ is Lyapunov for $X$ we 
immediately even get $\rho_x(\xi,X)\leq0$.
\item\label{Rema:6:ii}
In 
the case $\ind(x)=1$
we certainly have $\rho_x(\xi,X)\leq0$.
\item
Throughout the whole section we did not make use of a 
Morse--Smale condition.
\end{enumerate}
\end{remark}

Using Proposition~\ref{P:expgrow}, Proposition~\ref{rhoxprop}, 
Remark~\ref{rem_on_K}\itemref{rem_on_K:ii} and 
Remark~\ref{Rema:6}\itemref{Rema:6:ii} we get

\begin{corollary}
Suppose $\xi$ is Lyapunov for $X$ and $x$ a zero of $X$.
If $1<\ind(x)<\dim(M)$ assume moreover that $K_x(y)=\emptyset$ 
for all zeros $y$ of $X$ with $\dim(M)>\ind(y)\geq\ind(x)$. 
Then $X$ has the exponential growth property at $x$
and $\rho_x(\xi,X)<\infty$.
\end{corollary}
In particular if $\dim(M)=2$ the exponential growth conjecture is true.


\section{Proof of Theorems~\ref{T:2} and \ref{T:3}
}
  \label{S:thm23}

Let $X$ be a vector field with 
$\xi\in H^1(M;\mathbb R)$ a Lyapunov cohomology class.
Recall that in Section~\ref{S:intro} we have defined the instanton 
counting function (or the Novikov incidence)
$\mathbb I^{X,\mathcal O,\xi}_{x,y}:
\hat{\mathcal P}_{x,y}^\xi\to \mathbb Z$, cf.\ \eqref{E:5}.

The following proposition is a reformulation of a basic observation made
by S.P.~Novikov \cite{N93} in order to define his celebrated complex.

\begin{proposition}\label{P:8}
\
\begin{enumerate}
\item\label{P:8i}
For any $x\in\mathcal X_q$, $y\in\mathcal X_{q-1}$ and 
every real number $R$ the set
$$
\bigl\{\hat\alpha\in\hat{\mathcal{P}}_{x,y}^\xi\bigm| 
\mathbb I^{X,\mathcal O,\xi}_{x,y}(\hat{\alpha})\neq 0,
\overline\omega(\hat\alpha)\geq R\bigr\}
$$
is finite. Here $\omega$ is any closed one form representing $\xi$.
\item\label{P:8ii}
For any $x\in\mathcal X_q$, $z\in\mathcal X_{q-2}$ and
$\hat{\gamma}\in\hat{\mathcal P}_{x,z}^\xi$
one has
\begin{equation}\label{E:26}
\sum
\mathbb I^{X,\mathcal O,\xi}_{x,y}(\hat\alpha)\cdot
\mathbb I^{X,\mathcal O,\xi}_{y,z}(\hat\beta)=0.
\end{equation}
where the sum is over all $y\in\mathcal X_{q-1}$,
$\hat\alpha\in\hat{\mathcal P}_{x,y}^\xi$ and all
$\hat\beta\in\hat{\mathcal P}_{y,z}^\xi$ with
$\hat\alpha\star\hat\beta=\hat\gamma$.
\end{enumerate}
Formula~\eqref{E:26} implicitly states that the left side of the equality 
contains only finitely many non-zero terms.
\end{proposition}

Proposition~\ref{P:8} above is equivalent to 
Theorem~2 parts 1 and 2 in \cite{BH01}. 
The proof, originally 
due to Novikov can be also found in \cite{BH01}.

The following proposition will be the main
tool in the proof of Theorem~\ref{T:2}. 

\begin{proposition}\label{P:9}
Suppose $t\in\mathbb R$, 
$\omega$ a closed one form representing $\xi$ and
$t>\rho(\xi,X)$.
Then:
\begin{enumerate}
\item\label{P:9i}
For every $a\in\Omega^q(M)$ and every $x\in {\mathcal X}_q$ the integral
$$
\Int^q_{X,\omega, \mathcal O}(t)(a)(x)
:=\int_{W^-_x}e^{th_{x}}(i_x^-)^*a
$$
converges absolutely. In particular
it defines a
linear map $\Int^q_{X,\omega,\mathcal O}(t):
\Omega^q(M)\to\Maps(\mathcal X_q,\mathbb R)$.
\item\label{P:9ii}
The map $\Int^q_{X,\omega,\mathcal O}(t):\Omega^q(M)\to\Maps(\mathcal X_q,\mathbb R)$ 
is surjective
and
\begin{equation}\label{F:1}
\Int^{q+1}_{X,\omega,\mathcal O}(t)(d_\omega(t)a)(x)
=\sum_{y\in\mathcal X_q,\ \hat{\alpha}\in\hat{\mathcal P}_{x,y}^\xi}
e^{t\overline{\omega}(\hat{\alpha})}
\mathbb I^{X,\mathcal O,\xi}_{x,y}(\hat{\alpha}) 
\Int^q_{X,\omega,\mathcal O}(t)(a)(y)
\end{equation}
where the right side of \eqref{F:1} is a potentially infinite sum which is 
convergent.
\end{enumerate}
\end{proposition}

Recall that for an oriented $n$--dimensional manifold $N$ and
$a\in\Omega^n(N)$ one has
$|a|:=|a'|\Vol\in\Omega^n(M)$,
where $\Vol\in\Omega^n(N)$ is any volume form and
$a'\in C^\infty(N,\mathbb R)$ is the
unique function satisfying $a=a'\cdot\Vol$.
The integral $\int_N a$ is called absolutely convergent, if
$\int_N|a|$ converges.

The proof of Proposition~\ref{P:9} is given in details
in \cite{BH01} section~5, (cf.\ Proposition 4) and 
uses in an essential way Theorem~\ref{T:5} and Stokes' theorem.
Particular care is necessary in view of the fact that $W^-_x$ are not 
compact. The integration has to be performed on a non compact manifold 
and Stokes' theorem applied to non-compact manifolds with corners.

The proof of Theorem~\ref{T:2} boils down to the verification of the 
following claims:

\begin{claim}\label{cl:1}
For any $t>\sup\{\rho(\xi, X),T\}$, $x\in\mathcal X_{q+1}$ and
$y\in\mathcal X_q$ the possibly infinite sum
$$
\sum_{\hat\alpha\in\hat{\mathcal P}_{x,y}^\xi}
\mathbb I^{X,\mathcal O,\xi}_{x,y}(\hat\alpha)
e^{t\overline{\omega}(\hat\alpha)}
$$ 
is convergent and the formula 
\begin{equation}\label{E:22}
\delta_{X,\omega,\mathcal O}^q(t)(E_y)
=\sum_{x\in\mathcal X_{q+1},\ \hat\alpha\in\hat{\mathcal P}_{x,y}^\xi}
\bigl(\mathbb I^{X,\mathcal O,\xi}_{x,y}(\hat\alpha)
e^{t\overline\omega(\hat\alpha)}\bigr)E_x
\end{equation} 
defines a linear map $\delta_{X,\omega,\mathcal O}^q(t):
C^q(X)=\Maps(\mathcal X_q,\mathbb R)\to
C^{q+1}(X)=\Maps(\mathcal X_{q+1},\mathbb R)$ 
which makes 
$(C^*(X),\delta^*_{X,\omega,\mathcal O}(t))$ a smooth (actually analytic)
family of cochain complexes of finite dimensional Euclidean spaces.
Recall that $\{E_x\}_{x\in\mathcal X}$ denotes the characteristic 
functions of $x\in\mathcal X$ and 
$\{E_x\}_{x\in\mathcal X}$ provide the canonical base of $C^*(X)$ which, 
implicitly equips each component  
$C^q(X)$ with a scalar product, the unique scalar product which makes 
this base orthonormal. 
Recall also that 
$\overline\omega:\hat{\mathcal P}_{x,y}^\xi\to\mathbb R$ 
was defined in section~\ref{S:intro} before Proposition~\ref{Prop:1} 
and makes sense even when $\omega$ is not a representative of $\xi$ but 
still, its pull back on $\tilde M$ is exact.
\end{claim}

\begin{claim}\label{cl:2}
The linear maps $\Int^q_{X,\omega,\mathcal O}(t)$ are surjective and define a 
morphism of cochain complexes.
\end{claim}

\begin{claim}\label{cl:3}
There exists $T$ (larger than $\rho(\omega,X)$) so that for $t>T$
the linear map $\Int^q_{X,\omega,\mathcal O}(t)$ when restricted to 
$\Omega^q_\sm(M)(t)$ is an isomorphism and actually an $O(1/t)$--isometry.
\end{claim}

Everything but the $O(1/t)$--isometry statement in Claim~\ref{cl:3} is a 
straight forward consequence of Theorem~\ref{T:1} and Proposition~\ref{P:9} 
above.
To check this part of Claim~\ref{cl:3} we have to go back to the proof
of Theorems~3 and 4 in \cite{BH01}, section~6.
We observe that if $t$ is large enough the restriction of 
$\Int^q_{X,\omega,\mathcal O}(t)$ to the subspace $H_1(t)\subset \Omega^q(M)$ 
defined in \cite{BH01}, section~4 page~172 (cf.\ Proof of Theorem~3)
is surjective and then by Lemma~7 in \cite{BH01} so is the restriction of  
$\Int^q_{X,\omega,\mathcal O}(t)$ to $\Omega^q_\sm(M)(t)$. This because  
$H_1(t)$ and $\Omega^q_\sm(M)(t)$ are, 
by Lemma~7 in \cite{BH01} section~\ref{S:thm4}, as close as we want
for $t$ large enough. Since the spaces $\Omega^q_\sm(M)(t)$ and 
$C^q(X)$ have the same finite dimension, by the 
surjectivity in Claim~\ref{cl:2},  
$\Int^q_{X,\omega,\mathcal O}(t)$ is an isomorphism 
and, as shown in \cite{BH01} section~4 page~172 an $O(1/t)$--isometry.
We take as the base $E^\mathcal O_x(t)$ the differential forms 
$E^\mathcal O_x(t)=\Int^q_{X,\omega,\mathcal O}(t)^{-1}(E_x)$.
This finishes the proof of Theorem~\ref{T:2}. Theorem~\ref{T:3} is a 
consequence of Theorem~\ref{T:2} and of Claim~\ref{cl:3}.

We conclude this section with the following remarks. 
Let $X$ be a vector field which has $\xi$ as Lyapunov cohomology class.  
Suppose $X$ satisfies MS and $\rho(\xi,X)<\infty$.
Let $\omega$ be a closed one form representing $\xi$.

For $t>\rho(\xi,X)$ the finite dimensional vector spaces 
$$
C^q(X):=\Maps({\mathcal X}_q,\mathbb R)
$$
and the linear maps 
$$
\delta^q_{X,\omega,\mathcal O}(t):\Maps({\mathcal X}_q,\mathbb R)
\to\Maps({\mathcal X}_{q+1},\mathbb R)
$$ 
defined by  
$$
\delta_{X,\omega,\mathcal O}^q(t)(E_x)
:=\sum_{y\in\mathcal X_{q+1},\ \hat{\alpha}\in{\mathcal P}_{y,x}^\xi} 
\mathbb I^{X,\mathcal O,\xi}_{y,x}(\hat{\alpha})
e^{t\overline\omega(\hat\alpha)}E_y
$$
give rise to a cochain complex of finite dimensional vector spaces
$$ 
\mathbb C^*(X,\omega,\mathcal O)(t)
:=\{C^q(X),\delta^q_{X,\omega,\mathcal O}(t)\},
$$ 
and to a morphism of such complexes:
$$
\Int^*_{X,\mathcal O,\omega}(t):(\Omega^*(M),d_\omega(t))\to
\mathbb C^*(X,\omega,\mathcal O)(t)
$$
One can show 
(implicit in Theorem~\ref{T:3})
that 
$\Int^*_{X,\mathcal O,\omega}(t)$ induces an isomorphism in 
cohomology. This fact will be used in this paper only when  
$X$ is the gradient of a smooth function in which case it is a 
strightforward consequence 
of deRham's theorem with local coefficients.

Let $\omega_1$ and $\omega_2$ be two closed one forms representing the 
same cohomology class $\xi$ and let $f:M\to\mathbb R$ be a smooth function
so that $\omega_1-\omega_2=df$. The collections of linear maps  
$$
m^q_{f}(t):\Omega^q(M)\to\Omega^q(M),\qquad
m^q_{f}(t)(a):=e^{tf}a,
$$
where $a\in\Omega^q(M)$, and 
$$
s^q_{f}(t):\mathbb C^q(X,\omega_1,\mathcal O)\to
\mathbb C^q(X,\omega_2,\mathcal O),\quad
s^q_{f}(t)(E_x):=e^{tf(x)}E_x,
$$
where $E_x\in\Maps(\mathcal X_q,\mathbb R)$
denotes the characteristic function of $x\in\mathcal X_q$,
define morphisms of cochain complexes making the diagram 
$$
\begin{CD}
\bigl(\Omega^*(M),d_{\omega_1}(t)\bigr)
@>m^*_f(t)>> 
\bigl(\Omega^*(M),d_{\omega_2}(t)\bigr)
\\
@V\Int^*_{X,\mathcal O,\omega_1}(t)VV  @VV\Int^*_{X,\mathcal O,\omega_2}(t)V
\\
\mathbb C^*(X,\mathcal O,\omega_1)(t)  @>s^*_f(t)>> \mathbb C^*(X,\mathcal O,\omega_2)(t)
\end{CD} 
$$
commutative for any $t>\rho(\xi,X)$.

Indeed because $h^1_x-h^2_x=(f-f(x))\cdot i^-_x$ is 
bounded,
$\int_{W_x^-}e^{th^2_{x}}(i_x^-)^*a$ is absolutely convergent iff 
$\int_{W_x^-}e^{th^1_{x}}(i_x^-)^*a$ is. Here $h^1_x$ is associated to $\omega_1$ 
and $h^2_x$ to $\omega_2$.

\section{The regularization $R(X,\omega,g)$}
  \label{S:reg}

In this section we discuss the numerical invariant $R(X,\omega,g)$ associated 
to a vector field $X$, a closed one form $\omega$ and a Riemannian metric $g$.
The invariant is defined by a possibly divergent integral but  regularizable
and is implicit in the work of \cite{BZ92}. More on this invariant is contained 
in \cite{BH03}.

In section~\ref{S:intro4} we have considered the 
Mathai--Quillen form 
$\Psi_g\in\Omega^{n-1}(TM\setminus M;\mathcal O_M)$
of an $n$--dimensional Riemannian manifold $(M,g)$. 
The Mathai--Quillen form (see~\cite{MQ86}) is actually 
associated to a pair $\tilde\nabla=(\nabla,\mu)$ consisting 
of a connection and a parallel Euclidean structure on a 
vector bundle $E\to M$. If $E$ is of rank $k$ it is a $k-1$
form 
$\Psi_{\tilde\nabla}\in\Omega^{k-1}(E\setminus M;\mathcal O_E)$
with values in the pull back of the orientation bundle 
$\mathcal O_E$ of $E$ to the total space of $E$.
Here $M$ is identified with the zero section in the bundle $E$.
If $g$ is a Riemannian metric
let $\tilde\nabla^g:=(\nabla^g,g)$ denote the Levi--Civita pair
associated to $g$ and write $\Psi_g:=\Psi_{\tilde\nabla^g}$.

The Mathai--Quillen form has the following properties:
\begin{enumerate}
\item\label{MQa}
For the Euler form 
$E_{\tilde\nabla}\in\Omega^k(M;\mathcal O_E)$ associated to
$\tilde\nabla$ we have
$d\Psi_{\tilde\nabla}=\pi^*E_{\tilde\nabla}$.
\item\label{MQb}
For two $\tilde\nabla^1$ and $\tilde\nabla^2$ we have
$\Psi_{\tilde\nabla^2}-\Psi_{\tilde\nabla^1}
=\pi^*\cs(\tilde\nabla^1,\tilde\nabla^2)$
modulo exact forms. Here
$\cs(\tilde\nabla^1,\tilde\nabla^2)
\in\Omega^{k-1}(M;\mathcal O_E)/d\Omega^{k-2}(M;\mathcal O_E)$
is the Chern--Simon invariant.
\item\label{MQd}
For every $x\in M$ the form $-\Psi_{\tilde\nabla}$ restricts 
to the standard generator of 
$H^{k-1}(E_x\setminus 0;\mathcal O_E)$, 
where $E_x$ denotes the fiber over $x\in M$. 
Note that the restriction of $-\Psi_{\tilde\nabla}$ 
is closed by \itemref{MQa}.
\item\label{MQc}
Suppose $E=TM$, $\tilde\nabla^g$ is the Levi--Civita pair, 
and suppose that on the open set $U$ we have coordinates 
$x^1,\dotsc,x^n$ in which the Riemannian metric $g|_U$ is 
given by $g_{ij}=\delta_{ij}$. Then, with respect
to the induced coordinates 
$x^1,\dotsc,x^n,\xi^1,\dotsc,\xi^n$ on $TU$, the 
form $\Psi_g$ is given by 
$$
\Psi_g=\frac{\Gamma(n/2)}{2\pi^{n/2}}  
\sum_i(-1)^i\frac{\xi^i}{\bigl(\sum_j(\xi^j)^2\bigr)^{n/2}}
d\xi^1\wedge\cdots\wedge\widehat{d\xi^i}\wedge\cdots\wedge d\xi^n,
$$
cf.\ \cite{MQ86}.
\end{enumerate}

Let $X$ be a vector field on $M$, i.e.\ a section of the 
tangent bundle $TM$. We suppose that it has only isolated 
zeros, that is its zero set $\mathcal X$ is a discrete subset
of $M$. The vector field defines an integer valued map
$\IND:\mathcal X\to\mathbb Z$, where $\IND(x)$ denotes the 
Hopf index of the vector field $X$ at the zero $x\in\mathcal X$.
This integer $\IND(x)$ is the degree of the 
map $(U,U\setminus x)\to(T_xM,T_xM\setminus 0)$ obtained by 
composing $X:U\to TU$ with the projection $p:TU\to T_xM$ 
induced by a local trivialization of the tangent bundle 
on a small disk $U\subseteq M$ centered at $x$.

Choose coordinates around $x$ so that we can speak of the 
disk $U_\epsilon$ with radius $\epsilon>0$ centered $x$. It is 
well known that we have:
\begin{equation}\label{E:IND:PSI}
\IND_x=-\lim_{\epsilon\to0}\int_{\partial U_\epsilon}X^*\Psi_g
\end{equation}
Indeed, by \itemref{MQb} we may assume that $g$ is flat on 
$U_\epsilon$. Thus $E_g=0$ and $\Psi_g$ is closed on 
$U_\epsilon$ by \itemref{MQa}. Using \itemref{MQd} we see that
$-\Psi_g$ gives the standard generator of 
$H^{n-1}(TU_\epsilon\setminus U_\epsilon;\mathcal O_{U_\epsilon})$ 
and thus certainly 
$\IND(x)=-\int_{\partial U_\epsilon}X^*\Psi_g$.

The vector field $X$ has its rest points (zeros) non-degenerate
and in particular isolated, if the map $X$ is transversal to 
the zero section in $TM$. In this case $\mathcal X$ is an 
oriented zero dimensional manifold, whose orientation is 
specified by $\IND(x)$. Moreover we have
$$
\IND(x)=\sign\det H\in\{\pm1\},
$$
where $H:T_xM\to T_xM$ denotes the Hessian. Particularly, if 
there exist coordinates $x^1,\dotsc,x^n$ centered at $x$ 
so that
\begin{equation}\label{R:4}
X=-\sum_{1\leq i\leq k} x^i\tfrac\partial{\partial x^i}+ 
\sum_{i>k} x^i\tfrac\partial{\partial x^i}
\end{equation}
we get $\IND(x)=(-1)^k$.

Let $X^1$ and $X^2$ be two vector fields and 
$\mathbb X:=\{X_s\}_{s\in [-1,1]}$ a smooth homotopy from
$X^1$ to $X^2$, i.e.\ $X_s=X^1$ for $s\leq-1+\epsilon$ and
$X_s=X^2$ for $s\geq 1-\epsilon$.
The homotopy is called non-degenerate if the map 
$\mathbb X:[-1,1]\times M\to TM$ defined by 
$\mathbb X(s,x):=X_s(x)$ is transversal to the zero section
of $TM$. In this case necessarily 
$X^1$ and $X^2$ are vector fields with non-degenerate zeros 
and so are all but finitely many $X_s$. Moreover all $X_s$ 
have isolated zeros with indexes in $\{0,1,-1\}$ and the 
zero set $\tilde{\mathcal X}$ of $\mathbb X$ 
is an oriented one dimensional smooth submanifold of 
$[-1,1]\times M$. Note that we have
$$
\partial\tilde{\mathcal X}=\sum_{y\in\mathcal X^2}\IND(y)y
-\sum_{x\in\mathcal X^1}\IND(x)x.
$$
If $\mathbb X'$ is a second homotopy joining $X^1$ with $X^2$ 
then $\tilde{\mathcal X'}-\tilde{\mathcal X}$ is the boundary 
of a smooth $2$--cycle. Indeed, if we choose a homotopy of 
homotopies joining $\mathbb X$ with $\mathbb X'$ which is 
transversal to the zero section, then its zero set will do 
the job.

Given a closed one form $\omega$ on $M$ denote by 
$$
I_{\mathbb X,\omega}:=\int _{\tilde{\mathcal X}}p_2^*\omega,
$$
where $p_2:\tilde{\mathcal X}\to M$ denotes the restriction 
of the projection $[-1,1]\times M\to M$. It follows from the\
previous paragraph that $I_{\mathbb X,\omega}$ does not
depend on the homotopy $\mathbb X$ --- only on $X^1$, $X^2$ 
and $\omega,$ and therefor will be denoted from now on by $I(X^1, X^2,\omega).$

\begin{remark}\label{Rem:6}
If there exists a simply connected open set $V\subset M$ 
so that $\mathcal X_s\subset V$ for all $s\in[-1,1]$ then one 
can calculate $I_{\mathbb X,\omega}$ as follows:
Choose a smooth function $f:V\to\mathbb R$ so that 
$\omega|_V=df$. Then 
$$
I_{\mathbb X,\omega}
=\sum_{y\in\mathcal X^2}\IND(y)f(y)
-\sum_{x\in\mathcal X^1}\IND(x)f(x).
$$
The proof of this equality is a straight forward application of 
Stokes' theorem.
\end{remark}

With these considerations we will describe now the 
\emph{regularization} referred to in Section~\ref{S:intro4}, 
cf.\ \eqref{E:10}. First note that for a non-vanishing vector 
field $X$, a closed one form $\omega$ and a Riemannian metric 
$g$ the quantity 
\begin{equation}\label{R:1}
R(X,\omega,g):=\int_M\omega\wedge X^*\Psi_g
\end{equation}
has the following two properties. 
$$
R(X,\omega+df,g)-R(X,\omega,g)=-\int_MfE_g
$$
for every smooth function $f$.
If $g^1$ and $g^2$ are two Riemannian metrics then
$$
R(X,\omega,g^2)-R(X,\omega,g^1)
=\int_M\omega\wedge\cs(g^1,g^2)
$$
where $\cs(g^1,g^2)=\cs(\tilde\nabla^{g^1},\tilde\nabla^{g^2})$.
This follows from properties~\itemref{MQa} and \itemref{MQb} of 
the Mathai--Quillen form.

If $X$ has zeros, then the form $\omega\wedge X^*\Psi_g$
is well defined on $M\setminus\mathcal X$ but the integral 
$\int_{M\setminus\mathcal X}\omega\wedge X^*\Psi_g$ might
be divergent unless $\omega$ is zero on a neighborhood of 
$\mathcal X$.

We will define below a regularization of the integral 
$\int_{M\setminus\mathcal X}\omega\wedge X^*\Psi_g$
which in case $\mathcal X=\emptyset$ is equal to the 
integral~\eqref{R:1}. For this purpose we choose a smooth 
function $f:M\to\mathbb R$ so that the closed $1$--form 
$\omega':=\omega-df$ vanishes on a neighborhood of 
$\mathcal X$, and put
\begin{equation}\label{R:6}
R(X,\omega,g ;f)
:=\int_{M\setminus\mathcal X}\omega'\wedge X^*\Psi_g
-\int_MfE_g
+\sum_{x\in\mathcal X}\IND(x)f(x)
\end{equation}

\begin{proposition}
The quantity $R(X,\omega,g;f)$ is independent of $f$.
\end{proposition}

Therefore $R(X,\omega,g;f)$ can be denoted by $R(X,\omega,g)$
and will be called the \emph{regularization} of
$\int_{M\setminus\mathcal X}\omega\wedge X^*\Psi_g$.

\begin{proof}
Suppose $f^1$ and $f^2$ are two functions such that 
$\omega^i:=\omega-df^i$ vanishes in a neighborhood $U$ of 
$\mathcal X$, $i=1,2$. For every $x\in\mathcal X$ we 
choose a chart and let $D_\epsilon(x)$ denote the 
$\epsilon$--disk around $x$. Put 
$D_\epsilon:=\bigcup_{x\in\mathcal X}D_\epsilon(x)$.

For $\epsilon$ sufficiently small $D_\epsilon\subseteq U$ and 
$f^2-f^1$ is constant on each $D_\epsilon(x)$. 
From~\eqref{R:6}, Stokes' theorem and
\eqref{E:IND:PSI} we conclude that 
\begin{align*}
R(X,\omega,g;f^2)&-R(X,\omega,g;f^1)
-\sum_{x\in\mathcal X}\IND(x)\bigl(f^2(x)-f^1(x)\bigr)
=
\\
&=
-\int_{M\setminus\mathcal X}d\bigl((f^2-f^1)\wedge X^*\Psi_g\bigr)
\\
&=
-\lim_{\epsilon\to 0}\int_{M\setminus D_\epsilon}
d\bigl((f^2-f^1)\wedge X^*\Psi_g\bigr)
\\
&=
\sum_{x\in\mathcal X}\bigl(f^2(x)-f^1(x)\bigr)
\lim_{\epsilon\to 0}\int_{\partial D_\epsilon(x)}X^*\Psi_g
\\
&=
-\sum_{x\in\mathcal X}\IND(x)\bigl(f^2(x)-f^1(x)\bigr)
\end{align*}
and thus $R(X,\omega,g;f^1)=R(X,\omega,g;f^2)$.
\end{proof}

\begin{proposition}\label{P:15}
Suppose that $\mathbb X$ is a non-degenerate homotopy 
from the vector field $X^1$ to $X^2$ and $\omega$ is a 
closed one form. Then 
\begin{equation}\label{EE:38}
R(X^2,\omega,g)-R(X^1,\omega,g)=I(X^1, X^2,\omega).
\end{equation}
\end{proposition}

%
%

\begin{proof}
We may assume that there exists a simply connected 
$V\subseteq M$ with $\mathcal X_s\subseteq V$ for all 
$s\in[-1,1]$. Indeed, since both sides of \eqref{EE:38} do 
not depend on the homotopy $\mathbb X$
we may first slightly change the homotopy and assume that no 
component of $\tilde{\mathcal X}$ lies in a 
single $\{s\}\times M$. Then we find $-1=t_0,\dotsc,t_k=1$ so 
that for every $0\leq i<k$ we find a simply connected 
$V_i\subseteq M$ such that $\mathcal X_s\subseteq V_i$ 
for all $s\in[t_i,t_{i+1}]$.

Assuming $V$ as above we choose a function $f$ so that
$\omega':=\omega-df$ vanishes on a neighborhood of every 
$\mathcal X_s$, i.e.\ $p_2^*\omega'$ vanishes on a 
neighborhood of $\tilde{\mathcal X}$. Here
$p_2:[-1,1]\times M\to M$ denotes the canonical projection.
Moreover let $\tilde p_2:[-1,1]\times TM\to TM$ denote the 
canonic projection and note that 
$p_2^*\omega'\wedge\mathbb X^*\tilde p_2^*\Psi_g$ is a
globally defined form on $[-1,1]\times TM$.
Using Stokes' theorem and Remark~\ref{Rem:6} we then get:
\begin{eqnarray*}
R(X^2,\omega,g)-R(X^1,\omega,g)
-I_{\mathbb X,\omega}
&=&
\int_{[-1,1]\times M}d\bigl(
p_2^*\omega'\wedge\mathbb X^*\tilde p_2^*\Psi_g\bigr)
\\
&=&
\int_{[-1,1]\times M}p_2^*(\omega'\wedge E_g)
\\
&=&
0
\end{eqnarray*}
For the second equality we used 
$d\mathbb X^*\tilde p_2^*\Psi_g=p_2^*E_g$. The integrand of the
last integral vanishes because of dimensional reasons.
\end{proof}

With little effort, using Stokes' theorem and the properties
of the Mathai--Quillen form, one can proof 
$$
R(X,\omega+df,g)-R(X,\omega,g)
=-\int_MfE_g+\sum_{x\in\mathcal X}\IND(x)f(x)
$$
for every smooth function $f$, and
$$
R(X,\omega,g^2)-R(X,\omega,g^2)
=\int_M\omega\wedge\cs(g^1,g^2)
$$
for any two Riemannian metrics $g^1$ and $g^2$.
Its also not difficult to generalize the regularization
to vector fields with isolated singularities, cf.\ \cite{BH03}.

\section{Proof of Theorem~\ref{T:4}}
  \label{S:thm4}

The proof of Theorem~\ref{T:4} presented here combines results of 
Hutchings, Pajitnov and others (cf.\ \cite{H02}, \cite{P02}) with results 
of Bismut--Zhang, cf.\ \cite{BZ92}, 
\cite{BH03} and \cite{BFK01}. A recollection of these results, 
additional notations and preliminaries 
are necessary. They will be collected in four preliminary subsections. 
These subsections will be followed by the fifth where Theorem~\ref{T:4} is 
proven.

Recall from \cite{BFK01} that
a generalized triangulation $\tau=(f,g)$ on a closed manifold $M$ is a 
pair consisting of a Morse function $f$ and a Riemannian metric $g$ so 
that $X=-\grad_g f$ satisfies MS.

\subsection{Homotopy between vector fields}\label{S:homotopy}

Let $\xi\in H^1(M;\mathbb R)$,
and $\pi:\tilde M\to M$ be a covering so that $\pi^*\xi=0$.

Recall that a smooth family of sections $\mathbb X:=\{X_s\}_{s\in [-1,1]}$,
of the tangent bundle will be called a homotopy from the vector field 
$X^1$ to the vector field $X^2$ if there exists $\epsilon>0$ so that 
$X_s=X^1$ for $s<-1+\epsilon$ and $X_s=X^2$ for $s>1-\epsilon$.

To a homotopy $\mathbb X:=\{X_s\}_{s\in [-1,1]}$ one associates the 
vector field $Y$ on the compact manifold with boundary (cf 
appendix to section~\ref{S:top} for definition) $N:=M\times[-1,1]$
defined by 
\begin{equation}\label{A:1}
Y(x,s):= X(x,s)+1/2(s^2-1)\frac{\partial}{\partial s}.
\end{equation}
With this notation we have 
the following.

\begin{proposition}\label{P:13}
If $\mathbb X$ is a homotopy between two vector fields $X^1$ and $X^2$
which both have $\xi$ as a Lyapunov cohomology class. Then the vector 
field $Y$ has $p^*\xi$ as a Lyapunov cohomology class, 
cf.\ Definition~\ref{D:grad_with_bound}, where $p:N=M\times[-1,1]\to M$ 
is the first factor projection.
\end{proposition}

\begin{proof}
Since $X^1$ and $X^2$ are both vector fields with $\xi$ as Lyapunov cohomology class
we can choose closed $1$--forms $\omega_i$ representing $\xi$ and 
Riemannian metrics $g_i$ on $M$ such that $X^i=-\grad_{g_i}\omega_i$, 
$i=1,2$. Choose an admissible Riemannian metric $g$ on $N$ inducing 
$g_i$ on the boundaries; for example take
$$
g:=(1-\lambda)p^*g_1+\lambda p^*g_2+ds^2,
$$
where $\lambda:[-1,1]\to\mathbb R$ is a non-negative smooth function 
satisfying
$$
\lambda(s)=
\begin{cases}
\ 0 & \text{for $s\leq-1+\epsilon$ and}
\\
\ 1 & \text{for $s\geq1-\epsilon$.}
\end{cases}
$$

Next choose a closed $1$--form
$\omega$ on $N$ which restricts to $p^*\omega_1$ on
$M\times[-1,-1+\epsilon]$ and which restricts to $p^*\omega_2$ on  
$M\times[1-\epsilon,1]$. 
This is possible since $\omega_1$ and $\omega_2$ define the
same cohomology class $\xi$ and  can be achieved in the following way. 
Choose a function $h$ on $M$ with $\omega_2-\omega_1=dh$ and set
$\omega:=p^*\omega_1+d(\lambda p^*h)$. Choose a function
$u:[-1,1]\to\mathbb R$, such that:
\begin{enumerate}
\item
$u(s)=-\frac12(s^2-1)$ for all $s\leq-1+\epsilon$ and all $s\geq 1-\epsilon$.
\item
$u(s)\geq
\left\{\frac{-\omega(Y)(x,s)}{\frac12(s^2-1)}\right\}$
for all $s\in[-1+\epsilon,1-\epsilon]$ and all $x\in M$.
\end{enumerate}
This is possible since 
$\left\{\frac{-\omega(Y)(x,s)}{\frac12(s^2-1)}\right\}\leq0$
for $s=-1+\epsilon$ and $s=1-\epsilon$.

Then $\tilde\omega:=\omega+u(s)ds$ represents the cohomology
class $p^*\xi$ and one can verify that $Y$ is 
a vector field which coincides with $-\grad_{\tilde g}\tilde\omega$ 
in a neighborhood of $\partial N$ and for which $\tilde\omega(Y)<0$ 
on $N\setminus\mathcal Y$.
\end{proof}

Let $\mathbb X$ be a homotopy between two vector fields
$X^1$ and $X^2$ which satisfies MS. Let $Y$ be the vector field defined in
\eqref{A:1}. With the notations from the appendix to section~\ref{S:top}
(the case of a compact manifold with boundary) we have
$$
\mathcal Y=\mathcal Y''=\mathcal Y''_-\sqcup\mathcal Y''_+
$$
with
$$
\mathcal Y''_-=\mathcal X^1\times\{-1\}
\quad\text{and}\quad
\mathcal Y''_+=\mathcal X^2\times\{1\}.
$$

\begin{definition}
The homotopy $\mathbb X$ is called MS if the vector field $Y$ is 
MS, i.e.\ $X^1$ and $X^2$ are MS and for any $y\in\mathcal Y''_-$ and 
$z\in\mathcal Y''_+$ the maps $i^+_y$ and $i^-_z$ are transversal.
The homotopy $\mathbb X$ has exponential growth if $Y$ has exponential 
growth.
\end{definition}

\begin{proposition}\label{P:14}
Let $X^1$ and $X^2$ be two vector fields which satisfy MS and
$\mathbb X$ a homotopy from $X^1$ to $X^2$. Then there exists a MS
homotopy $\mathbb X'$ from $X^1$ to $X^2$, arbitrarily close 
to $\mathbb X$ in the $C^1$--topology.
\end{proposition}

\begin{proof} First we modify the vector field $Y$ 
into $Y'$ by a small change in the $C^1$--topology, and only in the 
neighborhood of $M\times\{0\}$, in order to have the Morse--Smale
condition satisfied for any $y\in\mathcal Y''_-$
and $z\in\mathcal Y''_+$. This can be done using Proposition~\ref{Prop:2}. 
Unfortunately $Y'$ might not have the $I$--component equal to 
$1/2(s^2-1)\partial/\partial s$, it is nevertheless $C^1$--close, so  
by multiplication with a function which is 
$C^1$--close to $1$ and equal to $1$ on the complement of a small 
compact neighborhood of the locus where $Y$ and $Y'$ are not the same,
one obtains a vector field $Y''$ whose $I$--component is exactly 
$1/2(s^2-1)\partial/\partial s$. The $M-$component of $Y''$ defines the 
desired homotopy. By multiplying a vector field with a smooth 
positive function the stable and unstable sets do not change, and their  
transversality continues to hold. 
\end{proof}

In view of Theorem~\ref{T:6} for compact manifolds with boundary 
we have the following

\begin{remark}\label{R:5}
For any $y=(x,1)\in\mathcal Y''_+$ the $1$--corner of $\hat W^-_y$ 
is given by
\begin{equation*}
\partial_1(\hat W^-_{y})
=V_0\sqcup V_1\sqcup V_2
\end{equation*} 
where
\begin{eqnarray*}
V_0 & \simeq & W^-_x
\\
V_1 & \simeq &
\bigcup_{
  {v\in\mathcal Y''_+,\hat\alpha\in\hat{\mathcal P}_{y,v}}
  \atop
  {\ind(v)=\ind(y)-1}
} 
\mathcal T(y,v,\hat\alpha)\times (W^-_v\setminus \partial N)
\\
V_2 & \simeq &
\bigcup_{
  {u\in\mathcal Y''_-,\hat\alpha\in\hat{\mathcal P}_{y,u}}
  \atop
  {\ind(u)=\ind(y)-1}
}
\mathcal T(y,u,\hat\alpha)\times W^-_u
\end{eqnarray*}
\end{remark}

It is understood that $W^-_x$ represents the unstable manifold 
in $M=M\times\{1\}$ if $x\in\mathcal X^2$.

In view of \eqref{A:1} we introduce the invariant 
$\rho(\xi,\mathbb X)\in\mathbb R\cup\{\pm\infty\}$ for any homotopy 
$\mathbb X$ by defining
$$
\rho(\xi,\mathbb X):=\rho(p^*\xi,Y)
$$
Clearly $\rho(\xi,\mathbb X)\geq\rho(\xi,X^i)$ for $i=1,2$.

Suppose $\mathbb X$ is a MS homotopy from the MS vector field $X^1$ to
the MS vector field $X^2$.
For each $X^i$ choose the orientations $\mathcal O^i$, $i=1,2$. 
Observe that the set
$\hat{\mathcal P}_{x',x}$ identifies to
$\hat{\mathcal P}_{(x',1),(x',-1)}$.
The orientations $\mathcal O^1$ and $\mathcal O^2$ define
the orientations $\mathcal O$ for the unstable manifolds of the rest points 
of $Y$. For $x^1\in\mathcal X^1$, $x^2\in\mathcal X^2$ and 
$\hat\alpha\in\hat{\mathcal P}_{x^2,x^1}$ define the incidences
\begin{equation}\label{A:12}
\mathbb I^{\mathbb X,\mathcal O^2,\mathcal O^1}_{x^2,x^1}(\hat\alpha):=
\mathbb I^{Y,\mathcal O}\bigl((x^2,1),(x^1,-1)\bigr)(\hat\alpha).
\end{equation}

Suppose in addition that $\rho(\xi,\mathbb X)<\infty$.
For any $t>\rho(\xi,\mathbb X)$ and $\omega$ a closed one form 
representing $\xi$ define the linear maps
$$
u^q_\omega(t):=u^q_{\mathbb X,\mathcal O^1,\mathcal O^2,\omega}(t):
\Maps(\mathcal X^1_q,\mathbb R)\to\Maps(\mathcal X^2_q,\mathbb R)
$$
and the linear maps
$$
h^q_{\omega}(t):=h^q_{\mathbb X,\mathcal O^1,\mathcal O^2,\omega}(t):
\Omega^q(M)\to\Maps(\mathcal X^2_{q-1},\mathbb R)
$$
by
\begin{equation}\label{A:13}
u^q_{\omega}(t)(E_{x^1}):=
\sum_{
  {x^2\in\mathcal X^2}
  \atop
  {\hat\alpha\in\hat{\mathcal P}_{x^2,x^1}}
}
\mathbb I^{\mathbb X,\mathcal O^2,\mathcal O^1}_{x^2,x^1}(\hat\alpha)
e^{t\omega(\hat\alpha)}E_{x^2},
\quad x^1\in\mathcal X^1_q
\end{equation}
and
$$
(h^q_{\omega}(t)(a))(x^2)
=\int_{W^-_y}
e^{t F_y}(i^-_y)^*p^*a,
\quad\text{$x^2\in\mathcal X^2_{q-1}$ and $y=(x^2,1)$.}
$$
The right side of \eqref{A:13} is a convergent infinite sum
since it is a sub sum of the right hand side of \eqref{E:22} when 
applied to the vector field $Y$.


\begin{proposition}
Suppose $X^1, X^2$ are two MS vector fields having $\xi$ as a Lyapunov 
cohomology class and suppose $\mathbb X$ is a MS homotopy. 
Suppose that $\rho:=\rho(\xi,\mathbb X)=\rho(p^*\xi,Y)<\infty$ and 
$\omega$ is a closed one form with $p^*\omega$ exact; here 
$p:\tilde M\to M$ is the $\Gamma$--principal covering associated with 
$\Gamma$. Then for $t>\rho$ we have:
\begin{enumerate}
\item\label{st:1}
The collection of linear maps $\{u^q_{\omega}(t)\}$ defines a
morphism of cochain complexes:
$$
u^*_\omega(t):=u^*_{\mathbb X,\mathcal O^1,\mathcal O^2,\omega}(t):
\mathbb C^\ast(X^1,\mathcal O^1,\omega)(t)\to\mathbb 
C^\ast (X^2,\mathcal O^2,\omega)(t)
$$
\item\label{st:2}
The collection of linear maps $h^q_{\omega}(t)$ defines an algebraic 
homotopy between 
$\Int^*_{X^2,\mathcal O^2,\omega}(t)$ 
and 
$u^*_{\mathbb X,\mathcal O^1,\mathcal O^2,\omega}(t)\circ
\Int^*_{X^2,\mathcal O^2,\omega}(t)$.
Precisely, we have:
$$
h^{*+1}_\omega(t)\circ d^*_\omega(t)
+\delta^{*-1}_\omega(t)\circ h^*_{\omega}(t)
=
u^*_{\omega}(t)\circ\Int^*_{X^1,\mathcal O^1,\omega}(t)
-\Int^*_{X^2,\mathcal O^2,\omega}(t)
$$
\end{enumerate}
\end{proposition}

\begin{proof}
Statement \itemref{st:1} follows from the equality
$$
\sum_{
  {x'\in\mathcal X^1_{q+1},\ \hat\alpha\in\hat{\mathcal P}_{z,x'}} 
  \atop
  {\hat\beta\in\hat{\mathcal P}_{x',x},\ \hat\alpha\star\hat\beta=\hat\gamma}
}
\mathbb I^{\mathbb X,\mathcal O^2,\mathcal O^1}_{z,x'}(\hat\alpha)
\mathbb I^{X^1,\mathcal O^1}_{x',x}(\hat\beta)
-
\sum_{
  \stackrel
  {z'\in\mathcal X^2_q,\ \hat\alpha\in\hat{\mathcal P}_{z,z'}}
  {\hat\beta\in\hat{\mathcal P}_{z',x},\ \hat\alpha\star\hat\beta=\hat\gamma}
}
\mathbb I^{X^2,\mathcal O^2}_{z,z'}(\hat\alpha)
\mathbb I^{\mathbb X,\mathcal O^2,\mathcal O^1}_{z',x}(\hat\beta) 
=0
$$
for any $x\in\mathcal X^1_q$, $z\in\mathcal X^2_{q+1}$ and
$\hat\gamma\in\hat{\mathcal P}_{z,x}$ which is a reinterpretation of 
equation \eqref{E:26} when applied to the vector field $Y$,
the rest points $(x,-1)$ and $(z,1)$ and 
$\hat\gamma\in\hat{\mathcal P}_{z,x}
=\hat{\mathcal P}_{(z,1),(x,-1)}$.
The sign stems from the fact that the sign associated to a trajectory
from $z$ to $z'$ changes when it is considered as trajectory in
$M\times[-1,1]$ instead of $M\times\{1\}$.

To verify statement \itemref{st:2} we first observe that:
\begin{enumerate}
\item[(a)]
If $y,u\in\mathcal Y$ the restriction of $F_y$ to
$\mathcal T(y, u)(\hat\alpha)\times W^-_u$, 
$\hat\alpha\in\hat{\mathcal P}_{y,u}$, 
when this set is viewed as 
a subset of $\hat W^-_y$ is given by 
$$
F_u\circ\pr_{W^-_u}+\overline\omega(\hat\alpha).
$$
\item[(b)]
If $y=(x,-1)$, $x\in\mathcal X^1$, via the identification of $W^-_x$ to 
$W^-_y$, we have $F_y=h_x$.
\end{enumerate}

In view of the uniform convergence of all integrals
which appear in the formulae below, guaranteed by the hypothesis $t>\rho$,
the Stokes theorem for manifolds with corners gives for any 
$a\in\Omega^q(M)$ and $y\in(\mathcal Y''_+)_q$
\begin{equation}\label{A:17}
\int_{\hat W^-_y}d(e^{tF_y}c)
=\int_{V_0}e^{tF_y}c 
+\int_{V_1}e^{tF_y}c
+\int_{V_2}e^{tF_y}c,
\end{equation}
where $c:=(i^-_{y})^*p^*a\in\Omega^q(\hat W^-_y)$.


In view of the Remark~\ref{R:5} we have
\begin{equation}\label{A:18}
\int_{V_0}e^{tF_y}c
=
\Int^q_{X^2,\mathcal O^2,\omega}(t)(a),
\end{equation}
and
\begin{eqnarray}
\notag
\int_{V_2}e^{tF_y}c
&=&
\sum_{
  {u\in{\mathcal Y}''_+,\ \hat\alpha\in\hat{\mathcal P}_{y,u}}
  \atop
  {\ind(u)=\ind(y)-1}
}
\mathbb I^{Y,\mathcal O}_{y,u}(\hat\alpha)e^{t\omega(\hat\alpha)}
\int_{\hat W^-_u}e^{tF_u}(\hat i^-_u)^*p^*a 
\\ \label{A:19}
&=&
-(\delta^{q-1}_{\omega}(t)\circ h^q_{\omega}(t))(a)
\end{eqnarray}
and
\begin{eqnarray}
\notag
\int_{V_1}e^{tF_{y}}c
&=&
-
\sum_{
  {v\in{\mathcal Y}''_-,\ \hat\alpha\in\hat{\mathcal P}_{y,v}}
  \atop
  {\ind(v)=\ind(y)-1}
} 
\mathbb I^{Y,\mathcal O}_{y,v}(\hat\alpha)e^{t\omega(\hat\alpha)}
\int_{\hat W^-_v}e^{tF_v}(i^-_v)^*p^*a
\\ \label{A:20}
&=&
-(u^q_\omega(t)\circ\Int^q_{X^1,\mathcal O^1,\omega}(t))(a).
\end{eqnarray}
Moreover
\begin{eqnarray}\label{A:21}
\notag
(h^{q+1}_{\omega}(t)\circ d_{\omega}(t))(b)(y)
&=&
\int_{\hat W^-_y}e^{tF_y}(i^-_y)^* p^*(db +t\omega \wedge b)
\\
&=&
\int_{\hat W^-_y}d(e^{tF_y}(i^-_y)^* p^* b) 
\end{eqnarray}
and the statement follows combining the 
equalities~\eqref{A:17}--\eqref{A:21}.
\end{proof}

The following proposition will be important in the proof of 
Theorem~\ref{T:4}.

\begin{proposition}\label{P:16}
\
\begin{enumerate}
\item\label{P:13:st:1}
Let $(f,g)$ be a pair consisting of a Morse function 
and a Riemannian metric. 
Then the vector field $-\grad_gf$ has any 
$\xi\in H^1(M;\mathbb R)$ as a Lyapunov cohomology class.
\item\label{P:13:st:2}
Let $X$ be a vector field which has MS property and has $\xi$ as Lyapunov 
cohomology class. Let
$\tau=(f,g)$ be a generalized triangulation. 
Then there exists a homotopy $\mathbb X$ from $X^1:=X$ to $X^2$ 
which is MS and is $C^0$--close to the family 
$X_s=\frac{1-s}2X-\frac{1+s}{2}\grad_gf$. One can choose $\mathbb X$ 
to be $C^1$--close to a family $l(s)X-(1-l(s))\grad_gf$ where 
$l:[-1,1]\to[0,1]$ is a smooth function with $l'(s)\leq0$ and $l'(s)=0$ 
in a neighborhood of $\{\pm 1\}$.
\end{enumerate}
\end{proposition}

\begin{proof}[Proof of \itemref{P:13:st:1}]

Let $\omega$ be a closed one form representing 
$\xi$ with support disjoint from a neighborhood of $\Cr(f)$.
Clearly for $C$ a large constant the form $\omega':=\omega+Cdf$ 
represents $\xi$ and satisfies $\omega'(-\grad_gf)<0$.   
\end{proof}

\begin{proof}[Proof of \itemref{P:13:st:2}]
First consider the family 
$X_s:=(\frac{1-s}2)X-(\frac{1+s}2)\grad_gf$. 
Change the parametrization to make this family 
locally constant near $\{\pm 1\}$, hence get a homotopy and apply 
Proposition~\ref{P:14} to change this homotopy into one which 
satisfies MS.
\end{proof}

\begin{definition}\label{D:strong_exp}
A vector field $X$ which satisfies MS and has $\xi$ as a 
Lyapunov cohomology class is said to have 
\emph{strong exponential growth} if 
for one (and then any) generalized triangulation $\tau=(f,g)$
there exists a homotopy $\mathbb X$ from $X$ to $-\grad_gf$ which has 
exponential growth, cf Definition 8, equivalently, with $\rho (\xi, \mathbb X),\infty.$ 
\end{definition}

To summarize the discussion in this subsection consider: 
\begin {enumerate}
\item 
a vector field $X^1=-\grad_{g'}\omega$ with $\omega$ a 
Morse form representing $\xi$ and $g'$ a Riemannian metric so 
that $X^1$ satisfies MS,  
\item 
A vector field $X^2=-\grad_{g''}f$, $\tau=(f,g'')$ a generalized 
triangulation,
\item 
A homotopy $\mathbb X$ from $X^1$ to $X^2$ which satisfies MS.
\end{enumerate}
Since $\rho(\xi,X^2)=-\infty$, for any $t\in\mathbb R$ we have a 
well defined morphism of cochain complexes 
$$
\Int^*_{X^2,\mathcal O^2,\omega}(t):(\Omega^*(M),d_\omega(t)) 
\to\mathbb C^\ast(X^2,\mathcal O^2,\omega)(t).
$$
If $X^1$ has $\rho(\xi,X^1)<\infty$, equivalently $X^1$ has 
exponential growth, then for $t$ large enough we have a 
well defined morphism of cochain complexes 
$$
\Int^*_{X^1,\mathcal O^1,\omega}(t):(\Omega^*(M),d_\omega (t)) 
\to\mathbb C^\ast(X^1,\mathcal O^1,\omega)(t). 
$$
If $\mathbb X$ has $\rho(\xi,\mathbb X)<\infty$ then for $t$ large enough 
we have the morphism of cochain complexes 
$$
u^*_{\mathbb X,\mathcal O^1,\mathcal O^2,\omega}(t):
\mathbb C^\ast(X^1,\mathcal O^1,\omega)(t)\to
\mathbb C^\ast(X^2,\mathcal O^2,\omega)(t)
$$  
and the algebraic homotopy $h^*_{\mathbb X,\mathcal O^1,\mathcal O^2,\omega}(t)$ 
making the diagram 
below homotopy commutative   
$$
\begin{CD}
(\Omega^*(M),d_{\omega}(t))
@>\Id>> 
(\Omega^*(M),d_{\omega}(t))
\\
@V\Int^*_{X^1,\mathcal O^1,\omega}(t)VV 
@VV\Int^*_{X^2,\mathcal O^2,\omega}(t)V
\\
\mathbb C^\ast(X^1,\mathcal O_1,\omega)(t) 
@>u^*_{\mathbb X,\mathcal O^1,\mathcal O^2,\omega}(t)>> 
\mathbb C^\ast(X^2,\mathcal O^2,\omega)(t),
\end{CD}
$$
with all arrows inducing isomorphisms in cohomology.

\subsection{A few observations about torsion}\label{SS:obser}

Consider cochain complexes $(C^*,d^*)$ of free $A$--modules of finite 
rank. The cohomology $H^*:= H^*(C^*,d^*)$ is 
also a graded $A$--free module
of finite rank. Here $A$ is a commutative ring with unit.

Recall that the bases $\underline m'\equiv\{m'_1,\dotsc,m'_k\}$ 
and $\underline m" \equiv\{m''_1,\dotsc,m''_k\}$ of the free $A$--module $M$ 
are equivalent iff the isomorphism
$T:M\to M$ defined by $T(m'_i)=m''_i$ has determinant $\pm 1$.

For two equivalence classes of bases, $[\underline c]$ of $C^*$ and $[\underline h]$ of $H^*,$
Milnor, cf.~\cite{Mi68}, has defined the torsion 
$\tau((C^*,d^*),[\underline c],[\underline h])\in A^+/\{\pm 1\}$ where $A^+$ denotes the 
multiplicative group of invertible elements of $A$.

Recall that the bases $\underline m'\equiv\{m'_1,\dotsc,m'_k\}$ 
and $\underline m" \equiv\{m''_1,\dotsc,m''_k\}$ of the free $A$--module $M$ 
are equivalent iff the isomorphism

If the complex $(C^*, d^*)$ is acyclic there is no need of the base $\underline h$ and one has 
$\tau((C^*,d^*),[\underline c])\in A^+/\{\pm 1\}.$ If $\alpha:A\to B$ is a unit 
preserving ring homomorphism, by tensoring $(C^*,d^*),[\underline c],[\underline h]$ with $B,$ 
regarded as an $A-$module
via $\alpha,$ one obtains $((C')^*,(d')^*),[\underline c'],[\underline h']$ a cochain 
complex of free $B-$modules whose cohomology 
is a free $B-$module and the bases $[\underline c'],[\underline h']$. Clearly
$$\tau(((C')^*,(d')^*),[\underline c'],[\underline h'])= 
\alpha(\tau((C^*,d^*),[\underline c],[\underline h])).$$

If $A$ is the field $\mathbb R$ or $\mathbb C$, hence $(C^*,d^*)$ is a 
cochain complex of finite dimensional vector spaces, 
and $\scalar$ are scalar products in $C^*$ one can define
the $T$--torsion, 
$T((C^*,d^*),\scalar)\in\mathbb R_+$, by the formula
\begin{equation*}
\log T((C^*,d^*),\scalar)=1/2\sum_i(-1)^{i+1}i\log\det{}'\Delta_i 
\end{equation*}
where $\det{}'\Delta_i$ is the product of the non-zero eigen values of 
$\Delta_i:=(d^{i+1})^\sharp\cdot d^i+d^{i-1}\cdot(d^{i})^\sharp$. 
Here $(d^{i})^\sharp$ denotes the adjoint of $d^i$ with respect of 
the scalar product $\scalar.$

If in addition a scalar product $\sscalar$ in cohomology $H^*:=H^*(C^*,d^*)$
is  given one defines  
the positive real numbers $\Vol(H^i,\scalar, \sscalar)$ to be the volume of 
the isomorphism 
$$\ker d^i /d^{i-1}(C^{i-1})\to H^i.$$
Here the first vector space is equipped with the scalar product induced 
from $\scalar$ and the second with the scalar product $\sscalar$.
(Recall that 
the volume of an isomorphism $\theta:(V_1,\scalar_1)\to(V_2,\scalar_2)$ 
between two Euclidean vector spaces is the positive real number defined by:
$\log\Vol(\theta):=1/2\log\det\theta^\sharp\cdot\theta.$)

If $A$ is $\mathbb R$ or $\mathbb C$ then any base $\underline c$ resp.\ $\underline h$
induce a scalar product $\scalar_{\underline c}$ resp.\ $\sscalar_{\underline h},$
the unique scalar product which makes the base orthonormal. Although 
equivalent bases do not necessary provide the same scalar products 
they do however lead to the same $T$--torsions. This follows by inspection from 
Milnor's definition. Moreover one has 
\begin{eqnarray*}
\lefteqn{|\tau((C^*,d^*),[\underline c],[\underline h])|=}
\\&=&
T((C^*,d^*),\scalar_{\underline c})
+\sum_i(-1)^i\log\Vol(H^i,\scalar_{\underline c}, \sscalar_{\underline h})
\end{eqnarray*}

Let $u^*:(C^*_1,d^*_1)\to(C^*_2,d^*_2)$ be a morphism of cochain 
complexes of free $A$--modules of finite rank which induce isomorphism 
in cohomology. Then the mapping cone $\mathcal Cu^*$ is an acyclic 
cochain complex of free $A$--modules of finite rank.

Two equivalence classes of bases $[\underline c_1]$ of $C^*_1$ and $[\underline c_2]$ 
of $C^*_2$ provide an equivalence class of bases $[\underline c]$ of 
$\mathcal Cu^*$, and permit to define 
$$
\tau(u^*,[\underline c_1],[\underline c_{2}]):=\tau(\mathcal Cu^*,[\underline c]).
$$
If $A$ is $\mathbb R$ or $\mathbb C$ the scalar products $\scalar_1$ and
$\scalar_2$ in $C^*_1$ and $C^*_2$ provide the scalar product $\scalar$ in 
$\mathcal Cu^*$ and permit to define  
$$
T(u^*,\scalar_1,\scalar_2):=T(\mathcal Cu^*,\scalar).
$$
If the scalar products $\scalar_i := \sscalar _{\underline c_i}$, $i=1,2$ 
are induced from the bases 
$c_i$, $i=1,2$ we also have 
\begin{equation}\label {A:41}
|\tau((u^*,[\underline c_1],[\underline c_2])|=T(u^*,\scalar_1,\scalar_2)
\end{equation}

It is a simple exercise of linear algebra (cf.~\cite{BFK01}) to check 
that:

\begin{proposition}\label{P:17}
\
\begin{enumerate}
\item\label{P:17:i}
Let $u^*:(C^*_1,d^*_1)\to(C^*_2,d^*_2)$
be a morphism of cochain complexes of free $A-$modules whose cohomology 
modules  $H^\ast_i= H^*(C^*_i, d^*_i),\ i=1,2 $
are  also free. 
Suppose that $u^\ast$ induces isomorphisms in cohomology and denote 
these isomorphisms by $H^*u.$
Suppose  
$\underline c_1$ and $\underline c_2$ are bases of $C^*_1$ and $C^*_2$
and $\underline h_1$ and $\underline h_2$ are bases in 
$H^\ast_1$ and $H^\ast_2.$ 
If 
$$\prod(\det H^iu)^{(-1)^i}=\pm1,$$   
with $(\det H^i u)$ 
calculated with respect to the bases $\underline h_1$ and $\underline h_2,$ then
$$
\tau(u^*,[\underline c_1],[\underline c_2])=\tau(C^*_2,[\underline c_2],[\underline h_2])\cdot 
\tau(C^*_1,[\underline c_1],[\underline h_1])^{-1}.
$$
\item\label{P:17:ii}
Suppose $A$ is $\mathbb R$ or $\mathbb C$ and $\scalar_1$ and $\scalar_2$
are scalar products on $C^*_1$ and $C^*_2$. Then 
\begin{eqnarray}\label{E:35}
\lefteqn{\log T(u^*,\scalar_1,\scalar_2)=\log T((C^*_1,d^*_1),\scalar_1)}
\\&&\notag
-\log\tau((C^*_2,d^*_2),\scalar_2)
+\log\Vol H^*u
\end{eqnarray}
where $\log\Vol H^*u=\sum_i(-1)^i\log\Vol H^iu$ and $\Vol H^iu$ 
is calculated with respect to the scalar product induced from 
$\scalar_i$, $i=1,2$ in cohomology.
Moreover if $u^*$ is an isomorphism then  
$$
\log T(u^*,\scalar_1,\scalar_2 )=\sum_i(-1)^i\log\Vol u^i.
$$
\end{enumerate}
\end{proposition}

We conclude this subsection by recalling the following result of 
Bismut--Zhang, see \cite{BZ92} and \cite{BFK01}.

Suppose that $(M,g)$ is a closed Riemannian manifold $X=-\grad_{g'}f$ 
with $\tau=(f,g')$ a generalized triangulation, $\rho$ 
a representation of $\pi_1(M)$ and $\mu$ a Hermitian structure in the 
flat vector bundle associated with $\rho$. Consider 
$\Int^*:(\Omega^*(M,\rho),d_\rho)\to(C^*(\tau, \rho),\delta_{\mathcal O,\rho})$
and equip each of these complexes with a scalar product,
the first complex with the scalar product induced from the Riemannian 
metric $g$ and the Hermitian structure $\mu$ and the second with the 
scalar product $\scalar_{\mu,\tau}$ induced from the generalized 
triangulation $\tau$ and the Hermitian structure $\mu$, cf.~\cite{BFK01}.
In this notation $C^q(\tau, \rho)$ can be viewed as the vector 
space of sections above $\cal X_q$ of the vector bundle $E_\rho\to M$ 
equipped with the hermitian structure 
$\mu.$ This is a finite dimensional vector space with a scalar product.

Equip the cohomology of these cochain complexes with the induced 
scalar product. Denote by $H^*\Int$ the isomorphism induced in 
cohomology and write 
\begin{equation}
\log VH(\rho,\mu,g,\tau)=\sum(-1)^q\log\Vol(H^q\Int)
\end{equation}

Let $\omega(\mu)$ be the closed one form induced by $\mu$ 
as described 
in \cite{BZ92} and \cite{BFK01}.
To recall its definition denote by  $\det E_\rho\to M$ the determinant line 
bundle of $E_\rho\to M$ and equip it with the induced flat connection 
and hermitian structure. The form $\omega(\mu)$ is the logarithmic differential 
of the norm (with respect to the induced hermitian structure) of a parallel section in 
$\det E_\rho\to M$. Such sections exist only locally but their logarithmic differential 
is defined 
on the entire manifold $M$ and independent of the choosen section.   

We have the following result due to 
Bismut--Zhang, cf.~\cite{BFK01}.

\begin{theorem}\label{T:7}
With the hypothesis above we have 
\begin{eqnarray*}
\lefteqn{\log T_\an(M,g,\rho,\mu)=
\log T(C^*(\tau, \rho),\delta_{\mathcal O,\rho},\scalar_{\mu,\rho})+}
\\&&+
\log VH(\rho,\mu,g,\tau)+ R(X,\omega(\mu),g)
\end{eqnarray*}
\end{theorem}

\subsection{A summary of Hutchings--Pajitnov results}\label{SS:hutch}

We begin by recalling the results of Hutchings cf.\ \cite{H02}. 
  
Let $M$ be a compact smooth manifold, $m\in M$ a base point and 
$\xi\in H^1(M;\mathbb R)$. Recall from section~\ref{S:intro1} that $\xi$ 
defines the free abelian group $\Gamma$ and induces the injective 
homomorphism $\xi:\Gamma\to\mathbb R$. Denote by $\pi:\tilde M\to M$ 
the principal $\Gamma$--covering canonically associated with 
$x$ and $\pi_1(M,m)\to\Gamma$. To $\xi$ we associate  
\begin{enumerate}
\item    
the Novikov ring $\Lambda_\xi$ with coefficients in $\mathbb R$
which is actually a field,
\item    
the subring $\Lambda_{\xi,\rho}\subset\Lambda_\xi$, for any 
$\rho\in [0,\infty)$, cf.\ below,  
\item    
the multiplicative groups of invertible elements 
$\Lambda^+_\xi\subset\Lambda_\xi$
and $\Lambda^{+}_{\xi,\rho}\subset\Lambda_{\xi,\rho}$.
\end{enumerate}

The Novikov ring $\Lambda_\xi$ consists of functions $f:\Gamma\to\mathbb R$
which satisfy the property that for any real number $R\in\mathbb R$ the 
cardinality of the set 
$\{\gamma\in\Gamma|f(\gamma)\neq 0,\xi(\gamma)\leq R\}$ is finite.
The multiplication in $\Lambda_\xi$ is given by convolution, 
cf.\ \cite{BH01}. We have also shown both in Section~\ref{S:intro} and in 
more details in \cite{BH01} how to interpret the elements of 
$\Lambda_\xi$ as Dirichlet series.
In this context $\Lambda_{\xi,\rho}$ is the subring of $\Lambda_\xi$ 
consisting of those elements whose corresponding Dirichlet series have 
the abscissa of convergence smaller than $\rho$.

Note $\mathbb Z[\Gamma]\subset\Lambda_{\xi,\rho}\subset\Lambda_\xi$ and 
$H^*_\sing(M;\Lambda_\xi):=H^*_\sing(\tilde M;\mathbb Z)
\otimes_{\mathbb Z[\Gamma]}\Lambda_\xi$ 
is a finite dimensional vector space over the field $\Lambda_\xi$.
Let $\det H^*_\sing(M;\Lambda_\xi)$ denote the one dimensional vector 
space over $\Lambda_\xi$ defined by 
$$\det H^*_\sing(M;\Lambda_\xi)= \bigotimes_i(\Lambda^i(H^i_\sing(M;\Lambda_\xi))^{\epsilon(i)}$$
where $V^{\epsilon (i)}= V$ if $i$ is even and $V^{\epsilon (i)}$ is the 
dual of $V$ if $i$ is odd.
Let $X$ be a vector field which satisfies MS and has $\xi$ as a 
Lyapunov cohomology class and let $\tilde X$ be the pullback of $X$ 
on $\tilde M$. Choose $\mathcal O$ a collection of orientations 
for the unstable manifolds of the rest points of $X$ and therefore of 
the rest points of $\tilde X$.  
Denote by $(NC^q_{X,\xi},\partial^q_\mathcal O)$ the Novikov cochain 
complex of free $\Lambda_\xi$ modules (vector spaces since $\Lambda_\xi$
is a field) as defined in \cite{BH01} and by $H^*_X(M;\Lambda_\xi)$ 
its cohomology. There exists a canonical isomorphism 
$$HV^*_X:H^*_X(M;\Lambda_\xi)\to H^*_\sing(M;\Lambda_\xi)$$ 
described below.

\emph{The isomorphism $HV^*_X:$}
Note that in view of Proposition~\ref{P:14}, for any two vector fields 
$X^1$ and $X^2$ which satisfy MS and 
have $\xi$ as Lyapunov cohomology class, there exists 
homotopies $\mathbb X$ from $X^1$ to $X^2$ which 
satisfy MS. The incidences, 
$\mathbb I^{\mathcal O^2,\mathcal O^1}_{x^2, x^1}(\hat \alpha)$
defined in subsection 6.1 provide 
a morphism 
$$
u^*_{\mathbb X,\mathcal O^1,\mathcal O^2}
:(NC^\ast_{X^1, \xi},\partial^\ast_{\mathcal O^1}) 
\to(NC^\ast_{X^2, \xi},\partial^\ast_{\mathcal O^2})
$$
of cochain complexes which induces isomorphism between their cohomology.
This cohomology isomorphism is independent of the homotopy 
$\mathbb X$
and will be denoted by 
$$
H^*u(X^1,X^2):H^*_{X^1}(M;\Lambda_\xi)\to H^*_{X^2}(M;\Lambda_\xi).
$$
To show this one proceeds in the following way. One introduces 
the concept of homotopy between two homotopies $\mathbb X^1$  and $\mathbb X^2$ 
(both homotopies from the vector field $X^1$ to the vector field $X^2$) and prove  
that it induces an algebraic homotopy from the morphism induced by 
$\mathbb X^1$ and $\mathbb X^2$, hence the same isomorphism in cohomology.

For any three vector fields $X^i,\ i=1,2,3$ which satisfy MS 
and have $\xi$ as a Lyapunov cohomology class one has 
$$
H^*u(X^3,X^2)\cdot H^*u(X^2,X^1)= H^*u(X^3,X^1).
$$

Let $\tau=(f,g)$ be generalized triangulation, and 
$\xi\in H^1(M;\mathbb R)$. Let $X':= -\grad_g f.$ By Proposition~\ref{P:16} 
$X$ has $\xi$ as 
Lyapunov cohomology class. The Novikov complex 
$(NC^q_{X',\xi},\partial^q_\mathcal O)$ identifies to the  
geometric cochain complex associated to $\tilde\tau=(\tilde f,\tilde g)$,
the pull back of $(f,g)$ to $\tilde M$, tensored 
(over $\mathbb Z[\Gamma]$) by $\Lambda_\xi$. Recall that the geometric 
(or Morse) cochain complex associated to $(\tilde f,\tilde g)$ is a  
cochain complex of free $\mathbb Z[\Gamma]$ modules which can be regarded 
as 
a quotient of 
$C^*_\sing(\tilde M).$ When tensored with $\Lambda_\xi$ over 
$\mathbb Z[\Gamma]$ it calculates the same cohomology as 
$C^*_\sing(\tilde M)\otimes_{\mathbb Z[\Gamma]} \Lambda_\xi.$
We have therefore a well defined isomorphism from $H^\ast_{X'}(M;\Lambda_\xi)$
to $H^*_\sing(M;\Lambda_\xi).$

The composition of the this isomorphism with $H^*u (X,X')$ provides the 
canonical isomorphism 
$HV^*_X:H^*_X(M;\Lambda_\xi)\to H^*_\sing(M;\Lambda_\xi)$. 

\vskip .2in  
Denote by $E(M,m)$ the set of Euler structures based at $m\in M$ cf.\ 
\cite{B99} or \cite{BH03} for a definition, and let $e\in E(M,m)$.
Recall that in the presence of $X$ an Euler structure $e$ is 
represented by an Euler chain (cf.\ \cite{BH03}) which consists of a 
collection of paths $\alpha_x$ from $m$ to $x\in \mathcal X$.
Each such path provides a lift $\tilde x$ of $x$ (i.e.\ $\pi(\tilde x)=x$)
and therefore a base $\{\tilde E_x|x\in\mathcal X\}$ with $\tilde E_x$
the characteristic function of $\tilde x$ regarded as an element of 
$NC^q_{X,\xi,\rho}$, $q=\ind(x)$. Conversely, any lift 
$s:\mathcal X\to\widetilde X$, $s(x)=\tilde x$ defines an Euler chain 
and therefore together with $X$ an Euler structure $e$. The path 
$\alpha_x$ is the image by $\pi$ of a smooth path from $m$ to $\tilde x$.
Different Euler chains representing the same Euler structure might provide
nonequivalent bases. All theses bases will be named  
\emph{$e$--compatible} and denoted by $\underline e.$
Any lift $s$ which defines with $X$ the Euler 
structure $e$ will be also called \emph{$e$--compatible.} 

Choose an element $o_H\in\det H^*_\sing(M;\Lambda_\xi)\setminus 0$,
and consider bases $h^*$ in $H^*(NC^*_{X,\xi},\partial^q_{\mathcal O})$ 
which represent via the isomorphism $HV^*_X$ the element $o_H$. 
They all will be called $o_H$--compatible. Again the $o_H$--compatible bases
might not be equivalent, however
an inspection of Milnor definition of torsion \cite{Mi68} implies that 
the element 
$$\tau((NC^q_{X,\xi},\partial^q_\mathcal O),[\underline e],[\underline h])
\in\Lambda^+_{\xi}/\{\pm 1\}$$ 
as defined in section~\ref{SS:obser} for $\underline e$ resp.\ $\underline h$ 
$e$--compatible resp.\
$o_H$--compatible bases depends only on $e$ and $o_H;$ 
therefore denoted by $\tau_\xi(X,e,o_H).$ 

If $H^*_\sing(M;\Lambda_\xi)=0$ there is no need of $o_H$ and we have 
$\tau_\xi(X,e)\in\Lambda^+_{\xi}/\{\pm 1\}$.
  
Note that if $X$ has exponential growth, in view of Theorem~\ref{T:3}, 
the complex 
$(NC^q_{X,\xi},\partial^q_\mathcal O)$ contains, for $\rho$ large enough,  
a subcomplex 
of free $\Lambda_{\xi,\rho}$ modules $(NC^q_{X,\xi,\rho},\partial^q_\mathcal O),$
cf. \cite {BH01}
so that 
$$
(NC^q_{X,\xi,\rho},\partial^q_\mathcal O)
\otimes_{\Lambda_{\xi,\rho}}\Lambda_\xi= 
(NC^q_{X, \xi},\partial^q_\mathcal O).
$$ 
Moreover an $e$--compatible base will provide a base of 
of $\Lambda_{\xi,\rho}$--modules in this subcomplex.
If $H^*_\sing(M;\Lambda_\xi)=0$ 
then $\tau_\xi(X,e)\in\Lambda^+_{\xi,\rho}/\{\pm 1\}$ for $\rho$ 
large enough. More general,
if $H^*_\sing(M;\Lambda_{\xi,\rho})$ is a free 
$\Lambda_{\xi,\rho}$--module and 
$o_H\in \det H^*_\sing(M;\Lambda_{\xi,\rho}):= 
\bigotimes_i(\Lambda^i(H^i_X(M;\Lambda_{\xi, \rho}))^{\epsilon(i)}$ we will have 
$\tau_{\xi}(X,e,o_H)\in \Lambda^+_{\xi,\rho}/\{\pm 1\}$.

If the homotopy $\mathbb X$ has exponential growth then,
for $\rho$ big enough, we have 
$u^q_{\mathbb X,\mathcal O^1,\mathcal O^2}
(NC^q_{X^1,\xi,\rho},\partial^q_{\mathcal O^1}) 
\subset(NC^q_{X^2,\xi,\rho},\partial^q_{\mathcal O^2})$
and $\tau(u^\ast_{\mathbb X,\mathcal O^1,\mathcal O^2},[\underline e_1],[\underline e_2])$ 
which depends only on $X^1$, $X^2$ and
$e$, lies in $\Lambda^+_{\xi,\rho}/\{\pm 1\}$.

Note that for $t>\rho$ we denote by $\ev_t:\Lambda_{\xi,\rho}\to\mathbb R$ the ring
homomorphism which associates to each $f\in\Lambda_{\xi,\rho}$ 
interpreted as a Dirichlet series $f$, the value of the Laplace 
transform $L(f)$ at $t$, cf.\ section~\ref{S:intro6}. When applied to torsion it calculates 
the torsion of the corresponding complex tensored by $\mathbb R$.

Suppose now that $X$ is MS and satisfies also NCT. As noticed in \cite{H02},  
$\mathbb Z_X\in \Lambda_\xi$ and then $e^{\mathbb Z_X}\in\Lambda_\xi^+$.
The main result of Hutchings can be formulated as follows 

\begin{theorem}\label{T:8}
If $X^1$ and $X^2$ are two vector fields which are MS and NCT and have 
$\xi$ as a Lyapunov cohomology class then 
$$
e^{\mathbb Z_{X^1}}\cdot\tau_\xi(X^1,e,o_H)
=e^{\mathbb Z_{X^2}}\cdot\tau_\xi(X^2,e,o_H).
$$
\end{theorem}
  
The proof of this theorem is given in \cite {H02}. The author considers 
only the acyclic case (in which case $o_H$ is not needed).
The acyclicity hypothesis is used only to insure that 
the Milnor torsion (cf.\ \cite{Mi68}) $\tau_\xi(X,e)$ can be defined. This can 
be also defined 
in the non-acyclic case at the expense of the orientation
$o_H$. The orientation $o_H$ induces via $v_X^*$ an orientation 
in the cohomology of the Novikov complex associated to $X$ and together 
with the Euler structure $e$ a class of bases in the Novikov complex. 
From this moment on the arguments in \cite{H02} can be repeated word by word. 
 \vskip .1in

Let $\mathbb X$ be a homotopy from the vector field $X^1$ to the 
vector field $X^2$ which is MS and suppose that both vector fields have 
$\xi$ as a Lyapunov cohomology class. The incidences 
$\mathbb I^{\mathcal O^1,\mathcal O^2}_{\cdots}$, cf.\ \eqref{A:12}, 
induced from $\mathbb X$ and the orientation 
$\mathcal O= \mathcal O^1 \sqcup \mathcal O^2$ provide a morphism 
$u^*_{\mathbb X,\mathcal O^1,\mathcal O^2}:
(NC^*_{X^1},\partial_{\mathcal O_1}^*)
\to(NC^*_{X^2 },\partial_{\mathcal O_2}^*)$ 
which induces an isomorphism in cohomology as already indicated.

Choose bases $\underline e_i$ in each of the Novikov complexes 
$(NC^*_{X^i},\partial^*_{\mathcal O^i})$, $i=1,2$, which are 
$e$--compatible. By the same inspection of the Milnor definition of 
torsion one concludes that 
$\tau (u^*_{\mathbb X},[\underline e_1],[\underline e_2])\in\Lambda_\xi^+/\{\pm 1\}$ 
defined in section~\ref{SS:obser} depends only on $X^1$, $X^2$ and $e$.  
In view of
Proposition~\ref{P:17}\itemref{P:17:i} 
and of Theorem~\ref{T:8} one obtains 
  
\begin{proposition}\label{P:18}
Suppose $X^2$ and $X^1$ are two vector fields which satisfy MS and NCT 
and have $\xi$ as a Lyapunov cohomology class. Let  
$e$ be an Euler structure as above. Then 
$$
\tau (u^*_{\mathbb X,\mathcal O^1,\mathcal O^2},[\underline e_1],[\underline e_2])=  
e^{\mathbb Z_{X^2}}\cdot e^{-\mathbb Z_{X^1}}.
$$
\end{proposition}

As a consequence  $\tau(u^*_{\mathbb X},[\underline e_1],[\underline e_2])$ 
depends only on $X^1$, 
$X^2$ and 
then can be denoted by 
$\tau (X^1,X^2)$.

Suppose $X^2=-\grad_{g''}f$, $\tau=(f,g'')$ a generalized triangulation.
\begin{corollary}\label{UM4}
$\tau(X^1,X^2)=e^{\mathbb Z_{X^1}}$.
\end{corollary}
   
It is not hard to see that Hutchings theorem is equivalent to this 
corollary. In this form the result was also established by Pajitnov 
\cite{P02}.

Suppose $X$ is a vector field with $\xi$ a Lyapunov cohomology class 
which satisfies MS and in addition has exponential growth. The exponential
growth implies that any $e$--compatible base of $NC^q_{X,\xi}$ is 
actually a base of $NC^q_{X,\xi,\rho}$ for $\rho$ large enough.
For $t>\rho$ the $\mathbb R$--linear maps 
$\Ev^\omega_t:NC^q_{X^1,\xi,\rho}\to\Maps(\mathcal X_q,\mathbb R)$
defined by 
$$
\Ev^\omega_t(f):=\sum_{\tilde x\in\pi^{-1}(x)}
f(\tilde x)e^{-t\tilde h(\tilde x)}
$$
provide a morphism of cochain complexes 
$\Ev^\omega_t:(NC^*_{X^1,\xi,\rho},\partial^*_{\mathcal O})
\to\mathbb C^*(X,\mathcal O,\omega)(t)$ 
with $\Ev^\omega_t(\tilde E_x)=e^{-t\tilde h(\tilde x)}E_x$.
We regard $\mathbb C^*(X,\mathcal O,\omega)(t)$ equipped with the canonical 
base $\{E_x\}$ with $E_x$ the characteristic function of $x\in\mathcal X$.
The isomorphism $\Ev^\omega_t$ factors through  an isomorphism from 
$(NC^*_{X^1,\xi,\rho},\partial^*_{\mathcal O})\otimes_{\ev_t}\mathbb R$
to $\mathbb C^*(X,\mathcal O,\omega)(t).$
If the Novikov complex $(NC^*_{X,\xi,\rho},\partial^*_{\mathcal O})$
is acyclic so is $\mathbb C^*(X,\mathcal O,\omega)(t)$.

Let $s:\mathcal X\to\widetilde{\mathcal X}$ be a compatible lift and 
$\underline e$ the associated $e$--compatible base. A simple inspection of 
Milnor definition of torsion leads to
\begin{eqnarray}
\ev_t(\tau(X,\xi,e))
&=&\label{E:32}
\ev_t(\tau((NC^*_{X,\xi,\rho},\partial^*_{\mathcal O}),[\underline e])
\\&=&\notag
\tau(\mathbb C^*(X,\mathcal O,\omega)(t),[E_x])
\cdot e^{-t\sum_{x\in\mathcal X}\IND(x)\tilde h(\tilde x)}
\end{eqnarray}
  
   
Suppose now that $X^i$, $i=1,2$, are two vector fields which have 
$\xi$ as a Lyapunov cohomology class and $\mathbb X$ is a homotopy 
from $X^1$ to $X^2.$
Suppose in addition that $X^i$ and $\mathbb X$ satisfy MS 
and have exponential growth. Then we obtain the morphism of 
Novikov cochain complexes 
$u^*_{\mathbb X,\mathcal O^1,\mathcal O^2}:(NC^*_{X^1,\xi,\rho},
\partial^*_{\mathcal O^1})\to 
(NC^*_{X^1,\xi,\rho},\partial^*_{\mathcal O^2})$ 
which induces an isomorphism in cohomology. When tensored 
by $\mathbb R$ via $ev_t:\Lambda_{\xi,\rho}\to \mathbb R$ this morphism identifies to 
$u^*_{\mathbb X,\mathcal O^1,\mathcal O^2,\omega}(t).$

The Euler structure $e\in E(M,p)$ permits to choose $e-$compatible lifts of
the rest points and 
then $e$--compatible bases $\underline e_1$ and $\underline e_2$.  
The inspection of Milnor definition of torsion leads to 
\begin{eqnarray}\label{E:50} 
\lefteqn{
\ev_t(\tau(u^*_{\mathbb X},\mathcal O^1,\mathcal O^2),[\underline e_1],[\underline e_2])=
} 
\\&&\notag
=\tau(u^*_{\mathbb X,\mathcal O^1,\mathcal O^2,\omega}(t),[E_{x_1}],[E_{x_2}])
\cdot e^{-t(\sum_{x\in\mathcal X^1}\IND(x)\tilde h(\tilde x)
-\sum_{x\in\mathcal X^2}\IND(x)\tilde h(\tilde x))}
\end{eqnarray}
where $E_{x_1}$ resp.\ $E_{x_2}$ are the canonical base provided by the 
rest points of $X^1$ resp.\ $X^2$.

Note that in view of Proposition~\ref{P:15} (Additional property) for any 
$e$--compatible lifts of $\mathcal X^1$ and $\mathcal X^2$
we have:
\begin{equation}\label{E:51}
\sum_{x\in\mathcal X^1}\IND(x)\tilde h(\tilde x)
-\sum_{x\in\mathcal X^2}\IND(x)\tilde h(\tilde x)
=I(X^1,X^2,\omega).
\end{equation}
Different lifts which define the same the Euler structure keep the left side of 
\eqref {E:51}n unchanged.

\subsection{The geometry of closed one form}\label{SS:closed}

Suppose $M$ is a connected smooth manifold and $p\in M$ is a base point.
The homomorphism $[\omega]:H_1(M;\mathbb Z)\to \mathbb R$ induces the one 
dimensional representation 
$\rho=\rho_{[\omega]}:\pi_1(M,p)\to\GL_1(\mathbb R)$ defined by 
the composition $\pi_1(M,p)\to H_1(M;\mathbb Z)\xrightarrow{[\omega]}
\mathbb R\xrightarrow{\exp}\mathbb R_+=\GL_1(\mathbb R)$.
The representation $\rho$ provides a flat rank one vector bundle
$\xi_\rho:E_\rho\to M$ with the fiber above $p$ identified with 
$\mathbb R$. This bundle is the quotient of trivial flat bundle 
$\widetilde{M}\times\mathbb R\to\widetilde{M}$ by the group 
$\Gamma$ 
which acts diagonally on $\widetilde{M}\times\mathbb R$.
Here $\widetilde{M}$ denotes the principal $\Gamma$--covering 
associated with $[\omega]$ and constructed canonically 
with respect to $p$ (from the set of continuous paths originating from $p$).
Note that $\tilde M$ is equipped with a base point $\tilde p$ 
corresponding to the constant path in $p$. 
The group $\Gamma$ acts 
freely on $\widetilde M$ with quotient space $M$. The action of 
$\Gamma$ on $\mathbb R$ is given by the representation 
$\Gamma\to\mathbb R\xrightarrow{\exp}\mathbb R_+=\GL_1(\mathbb R)$.

There is a bijective correspondence between the closed one forms 
$\omega$ in the cohomology class represented by $\rho$ and the Hermitian
structures $\mu$ in the vector bundle $\xi_\rho$ which agree with a given
Hermitian structure on the fiber above $p$ (identified to $\mathbb R$).

Given $\omega$ in the cohomology class $[\omega]\in H^1(M;\mathbb R)$, 
one constructs a Hermitian structure $\tilde{\mu}_{\omega}$ on the trivial
bundle $\widetilde{M}\times\mathbb R\to\widetilde{M}$ which is 
$\Gamma$--invariant. Therefore, by passing to quotients one obtains a 
Hermitian structure $\mu_{\omega}$ in $\xi_\rho$. The Hermitian structure
$\tilde\mu_\omega$ is defined as follows: 
\begin{enumerate}
\item
Observe that the pull back $\tilde{\omega}$ of $\omega$ on 
$\widetilde{M}$ is exact and therefore equal $d\tilde h$ where 
$\tilde h:\widetilde M\to\mathbb R$ is the unique function with 
$\tilde h(\tilde {p})=0$ and $d\tilde h=\tilde\omega$.
\item
Define $\tilde{\mu}(\tilde x)$ by specifying the length of the vector 
$1\in\mathbb R$. We put $||1_{\tilde x}||_{\tilde\mu(\tilde x)}:=
e^{\tilde h(\tilde x)}$.
\end{enumerate}
Given a Hermitian structure $\mu$ one construct a closed one 
form $\omega_\mu$ as follows: Denote by 
$(\tilde E_\rho\to\tilde M,\tilde\mu)$ the pair consisting of the flat 
line bundle $\tilde E_\rho\to\tilde M$ and the Hermitian structure 
$\tilde\mu$ the pullback of the pair $(E_\rho\to M,\mu)$ 
to $\tilde M$ by the map $\widetilde M \to M$.
Let $\overline\mu$ be the Hermitian structure obtained by parallel 
transporting the scalar product $\tilde \mu_{\tilde p}$.
Denote by $\alpha:\tilde M\to\mathbb R$ the function 
$\alpha(\tilde x):=||v||_{\tilde\mu(\tilde x)}/||v||_{\overline{\mu}(\tilde x)}$ 
for $v$ a nonzero vector in $\tilde E_{\tilde x}$.

Define $\tilde\omega_{\mu}:=d\log(\alpha)$ and observe that this is a 
$\Gamma$ invariant closed one form, hence descends to a closed one form 
$\omega$ on $M$.

To simplify the writing below we denote (by a slight abuse of notation):  
\begin{enumerate}
\item 
$\rho(t):=\rho_{t\omega}$ 
\item 
$\mu(t):=\mu_{t\omega}$
\end{enumerate}

\begin{remark}\label{R:9}
The cochain complex $(\Omega^*(M),d_{\omega}(t))$ equipped with the 
scalar product induced from $g$ is isometric to the cochain complex 
$(\Omega^*(M,\rho(t))$ equipped with the scalar product induced from 
$g$ and $\mu(t)=\mu$ as defined in \cite{BFK01}.
\end{remark}

In particular we have 
$$
\log T^{\omega,g}_\an(\omega,t)=\log T_\an(M,\rho(t),g,\mu(t))
$$
where $\log T_\an(M,\rho,g,\mu)$ is the analytic torsion considered 
in \cite{BFK01} and associated with the Riemannian manifold $(M,g)$ 
the representation $\rho$ and the Hermitian structure $\mu$ in the flat 
bundle induced from $\rho$.

\begin{remark}\label{R:10}
Let $\xi\in H^1(M;\mathbb R)$ and $\omega$ a closed one form 
representing $\xi$. Suppose $X=-\grad_{g''}f$ where $\tau=(f,g'')$ is 
a generalized triangulation. Choose orientations $\mathcal O$ for the 
unstable manifolds of $X$. The morphism 
$$
\Int^*_{X,\mathcal O,\omega,}(t):(\Omega^*(M),d_\omega(t))
\to\mathbb C^*(X,\mathcal O,\omega)(t)
$$
defined in \eqref{E:9} where $(\Omega^*(M)$ is equipped with the scalar 
product induced from $g$ and  $\Maps(\mathcal X_q,\mathbb R)$ is equipped
with the obvious scalar product, i.e.\ associated with the base $\{E_x\}$,
is isometrically conjugate to 
$$
\Int^*:(\Omega^*(M,\rho(t)),d_{\rho(t)})
\to(C^*(\tau, \rho(t)),\delta_{\mathcal O,\rho(t)})
$$
defined in \cite{BFK01} where $(\Omega^{*}(M,\rho(t))$ is equipped with 
the scalar product induced by $(g,\mu(t))$ and $C^*(\tau, \rho(t))$ is equipped 
with the scalar product induced from $\tau$ and $\mu(t)$. 
\end{remark}

\subsection{Proof of Theorem~\ref{T:4}}

We begin with a triple $(g,g',\omega)$ with $X^1=X=-\grad_{g'}\omega$ as 
in the hypothesis of Theorem~\ref{T:4}. We choose orientations 
$\mathcal O^1$ for the unstable manifolds of $X^1$
We also choose $X^2=-\grad_{g''}f$ so that $\tau=(f,g'')$ is
a generalized triangulation and choose orientations $\mathcal O^2$ for 
the unstable manifolds of $X^2$.

For simplicity of the writing we will use the following abbreviations:
$I_1(t):=\Int^*_{X^1,\mathcal O^1,\omega}(t)|_{\Omega^{*}_\sm(M)}$ and
$I_2(t):=\Int^*_{X^2,\mathcal O^2,\omega}(t)|_{\Omega^{*}_\sm(M)}$.

In view of Proposition~\ref{P:17}\itemref{P:17:ii} applied to $I_1(t)$
one obtains 
\begin{eqnarray}\label{E;37}
\log(\mathbb V(t))
&=&
\log T(\mathbb C(X^1,\mathcal O^1,\omega)(t),\scalar_1)
\\&&\notag
-\log T^{\omega,g}_{\an,\sm}(t)
+\log\Vol(H^*(I_1(t))
\end{eqnarray}
where $\scalar_1$ is the scalar product induced from the canonical base 
$\{E_x\}$, $x\in\mathcal X^1$.

In view of Theorem~\ref{T:7} and the 
Remarks~\ref{R:9} and \ref{R:10} in section~\ref{SS:closed} one has 
\begin{eqnarray}\label{E;38}
\log T^{\omega,g}_\an(t)
&=&
\log T(\mathbb C(X^2,\mathcal O^2,\omega)(t),\scalar_2)
\\&&\notag
+\log\Vol(H^*(I_2(t))-tR(\omega,g,X^2) 
\end{eqnarray}
where $\scalar_2$ is the scalar product induced from the 
canonical base $\{E_x\}$, $x\in\mathcal X^2$.

Combining with \eqref{E;37} and \eqref{E;38} one obtains
\begin{eqnarray}\label{E;39}
\lefteqn{
\log(\mathbb V(t))-\log T^{\omega,g}_{\an,\la}(t)=}
\\
&=&\notag
\log\Vol(H^*(I_1(t))-\log\Vol(H^*(I_2(t))
\\&&+\notag
\log\tau(\mathbb C(X^1,\mathcal O^1,\omega)(t),\scalar_1)
\\&&-\notag
\log T(\mathbb C(X^2,\mathcal O^2,\omega)(t),\scalar_2)
+tR(X^2,\omega,g)
\end{eqnarray}

First consider the case that $X=X^1$ has exponential growth and 
$H^*(M,t[\omega])=0$ for $t$ large enough.
Note that $X^2$ 
has exponential growth too by Proposition~\ref{P:16}.  Clearly then 
$\log\Vol(H^*(I_1(t))=\log\Vol(H^*(I_2(t))=0$.
By \eqref{A:41}, \eqref{E:32} and \eqref{E:51} we have 
\begin{eqnarray}
\notag
\lefteqn{
\log T (\mathbb C(X^1,\mathcal O^1,\omega)(t),\scalar_1)
-\log T(\mathbb C(X^2,\mathcal O^2,\omega)(t),\scalar_2)=}
\\
&=&\notag
\log(\ev_t(\tau(X^1,\xi,e^*_1))
\\&&\label{E:44}
-\log(\ev_t(\tau(X^2,\xi,e^*_2))-
tI(X^1,X^2,\omega)
\end{eqnarray}
By Theorem~\ref{T:8} in section~\ref{SS:hutch} we have 
\begin{eqnarray}
\notag
\lefteqn{\log(\ev_t(\tau(X^1,\xi,e^*_1))
-\log(\ev_t(\tau (X^2,\xi,e^*_2))=}
\\
&=&\notag
\log(\ev_t(e^{\mathbb Z_{X^1}}\cdot e^{-\mathbb Z_{X^2}}))
\\&=&\notag
\log(\ev_te^{-\mathbb Z_{X^1}})
\\&=&\label{E:43}
-L(\mathbb Z_{X^1})(t)
\end{eqnarray}

Combining \eqref{E;39}, \eqref{E:44} and \eqref{E:43} one obtains 
the result.

Second consider the case $X$ has (strong) exponential growth property.
Then choose a homotopy $\mathbb X$ which satisfy:
$\rho(\xi,\mathbb X)<\infty$. 
Then for $t$ big enough,  
we have the following (algebraically) homotopy commutative diagram of 
finite dimensional cochain complexes 
whose arrows induce isomorphisms in cohomology. 
$$
\begin{CD}
(\Omega^*_\sm(M)(t),d_{\omega}(t))
@>\Id>>(\Omega^*_\sm(M)(t),d_{\omega}(t))
\\
@V\Int^*_{X^1,\mathcal O,\omega_1}(t)VV
@VV\Int^*_{X^2,\mathcal O,\omega_2}(t)V
\\
\mathbb C^\ast(X^1,\mathcal O_1,\omega)(t)  
@>u^*_{\mathbb X,\mathcal O_1,\mathcal O_2,\omega}(t)>> 
\mathbb C^\ast(X^2,\mathcal O_2,\omega)(t)
\end{CD} 
$$

For simplicity we write
$u^*(t):=u^*_{\mathbb X,\mathcal O^1,\mathcal O^2,\omega}(t)$
and observe that in view of the homotopy commutativity of the above 
diagram and of Proposition~\ref{P:17}\itemref{P:17:ii} we have 
\begin{eqnarray}\label{E:40}
\lefteqn{\log T(u^*(t),\scalar_1,\scalar_2)=}
\\&=&\notag
\log\Vol(H^*(I_1(t))-\log\Vol(H^*(I_2(t))
\\&&\notag
+\log T(\mathbb C(X^2,\mathcal O^1,\omega)(t),\scalar_1)
\\&&\notag
-\log T(\mathbb C(X^2,\mathcal O^1,\omega)(t),\scalar)
\end{eqnarray}

As noticed $(\mathbb C^\ast(X^2,\mathcal O^2,\omega)(t),\scalar)$ is 
isometric to $(\mathbb C(M,\tau,\rho(t),\mu(t))$.

By \eqref{A:41}, \eqref{E:50} and Proposition~\ref{P:18} combined with 
the observations that $X^2$ has no closed trajectories we have 
\begin{eqnarray}\label{E:41}
\log T(u^*(t),\scalar_1,\scalar_2)
&=&\notag
\log\tau(u^*(t),e^*_1,e^*_2)  
\\&=&
\ev_t(\mathbb Z_{X^1})
+ I(X^1,X^2,t\omega)
\end{eqnarray}
Combining \eqref{E;39} and \eqref{E:41} we obtain 
\begin{equation}\label{E:42}
\log\mathbb V(t)+\ev_t(\mathbb Z_{X^1})+ tI(X^1,X^2,\omega)
= 
\log T^{\omega, g}_{\an,\la}(t)+tR(X^2,\omega,g) 
\end{equation}
which in view of Proposition~\ref{P:15} implies the result.

When $H^*(M;t\Lambda_\xi)$ is acyclic we do not need the morphism 
$u^*$ and a simple consequence of Corollary~\ref{UM4} implies the result.
It turns out the \emph{strong exponential growth} can be weaken to the 
(apparently) weaker hypothesis $H^*_\sing(M;\Lambda_{\xi,\rho})$ is a 
free module over $\Lambda_{\xi,\rho}$ for some $\rho$.  
As in acyclic case one can circumvent the morphism $u^*(t)$.

%
%
%

\end{document}